\documentclass[11pt]{article}
\usepackage{amssymb,amsmath,amsfonts,latexsym,psfrag,graphicx,epsf} 
\newtheorem{lemma}{Lemma}[section]
\newtheorem{theorem}[lemma]{Theorem}

\newtheorem{corollary}[lemma]{Corollary}
\newtheorem{proposition}[lemma]{Proposition}

\newtheorem{definition}[lemma]{Definition}

\newtheorem{conjecture}[lemma]{Conjecture}
\newtheorem{problem}[lemma]{Problem}

\newcommand{\comment}[1]{}

\def\binom#1#2{{#1\choose#2}}
\def\ex{{\text{\rm ex}}}
\def\Grr{{{\mathbb G}_{r\times r}}}
\def\G33{{{\mathbb G}_{3\times 3}}}
\def\I2{{\mathbb I}_{\geq 2}}
\def\M5{{\mathbb M}_5}
\def\P4{{\mathbb P}_4}
\def\qed{\ifhmode\unskip\nobreak\hfill$\Box$\medskip\fi\ifmmode\eqno{\Box}\fi}

\def\cA{{\mathcal A}}

\def\cB{{\mathcal B}}

\def\bbE{{\mathbb E}}

\def\bbF{{\mathbb F}}

\def\cF{{\mathcal F}}

\def\bbG{{\mathbb G}}

\def\bH{{\mathbf  H}}
\def\cH{{\mathcal H}}

\def\cI{{\mathcal I}}

\def\cM{{\mathcal M}}
\def\bbP{{\mathbb P}}

\def\cP{{\mathcal P}}

\def\bR{{\mathbf  R}}
\def\cR{{\mathcal R}}

\def\bbT{{\mathbb T}}

\begin{document}

\pagestyle{myheadings}
\markright{{\small \sc Z. F\"uredi \& M. Ruszink\'o:}
  {\it\small Uniform Hypergraphs Containing no Grids}}
\thispagestyle{empty}

\title{Uniform Hypergraphs Containing no Grids
\footnote{ This copy was printed on {\today}.\quad
    {\rm\small {\jobname}.tex,} \hfill Version as of 
    Feb 24, 2011.
\break\indent{\it Keywords:} Tur\'an hypergraph problem, density problems, union
free hypergraphs, superimposed codes.
\hfill \break\indent{\it 2010 Mathematics Subject Classification:} 05D05, 11B25.} }
\author{{\bf Zolt\'an F\"uredi}
\thanks{ Research supported in part by the Hungarian National Science Foundation
 OTKA, and by the National Science Foundation under grant NFS DMS 09-01276.}
\\ Department of Mathematics, University of Illinois at Urbana-Champaign,
\\ Urbana, IL 61801, USA \quad and
\\ R\'enyi Institute of Mathematics of the Hungarian Academy of Sciences,
\\ Budapest, P. O. Box 127, Hungary-1364
\\ e-mail: {\tt z-furedi@illinois.edu}\quad  and\quad {\tt furedi@renyi.hu}
\\ and
\\ {\bf Mikl\'os Ruszink\'o}
\thanks{ Research supported in part by OTKA Grant K68322.}
\\ Computer and Automation Research Institute of the Hungarian
Academy of Sciences,
\\ Budapest, P. O. Box 63, Hungary-1518
\\ e-mail: {\tt ruszinko@sztaki.hu}
}

\date{${}$}
\maketitle
\begin{abstract}
A hypergraph is called an $r\times r$ {\em grid}  if it is isomorphic to a pattern of
  $r$ horizontal and $r$ vertical lines, i.e.,
 a family of sets   $\{A_1, \dots ,A_r, B_1,\dots ,B_r\}$ such that
 $A_i\cap A_j=B_i\cap B_j=\emptyset$ for $1\le i<j\le r$ and
 $|A_i\cap B_j|=1$ for $1\le i,j\le r$.
Three sets $C_1,C_2,C_3$
form a  {\em triangle} if they pairwise intersect in three distinct
singletons, $|C_1\cap C_2|=|C_2\cap C_3|=|C_3\cap C_1|=1$, $C_1\cap C_2\neq C_1\cap C_3$.
A hypergraph is {\em linear}, if $|E\cap F|\leq 1$ holds for every pair of edges.

In this paper we construct large linear $r$-hypergraphs which
contain no grids. Moreover, a similar construction gives large
linear $r$-hypergraphs which contain neither grids nor triangles.
For $r\ge 4$ our constructions are almost optimal. These
investigations are also motivated by coding theory: we get new
bounds for optimal superimposed codes and designs.
\end{abstract}

\section{Sparse hypergraphs, designs, and codes}
\label{intro} In this section we first present some previous
investigations in extremal set theory on the topic described in the
abstract. Then we state our main theorem. This is followed by
motivations in coding theory and corollaries where we improve the
previously known bounds for so called {\em optimal} superimposed
codes and designs. To prove the main theorem we are using tools from
combinatorial number theory and discrete geometry given in Section
2. In Section 3 we present constructions proving the stated
theorems, followed by remarks on union-free and cover-free graphs
and triple systems.

\subsection{Avoiding grids in linear hypergraphs}

Speaking about a  hypergraph $\mathbb F=(V, {\mathcal F})$ we frequently
 identify the vertex set $V=V(\bbF)$  by the set of first integers $[n]:=\{ 1,2, \dots, n\}$,
 or points on the plane $\bR^2$, or elements of a $q$-element finite field $F_q$.
To shorten notations we frequently say  'hypergraph $\cF$' (or set system $\cF$) thus
 identifying $\bbF$ to its edge set $\cF$.
$\bbF$ is {\em linear} if for all $A,B\in {\cal F}$, $A\neq B$ we have $|A\cap B|\le 1$.
The {\em degree}, $\deg_\bbF(x)$,  of an element $x\in [n]$ is the number of hyper-edges in ${\cal F}$ containing $x$.
${\bbF}$ is {\em regular} if every element $x\in [n]$ has the same degree.
It is {\em uniform} if every edge has the same number of elements,
 $r$-uniform means $|F|=r$ for all $F\in \cF$.
An $(n,r,2)$-{\em packing}  is a linear $r$-uniform hypergraph $\cP$ on $n$ vertices.
Obviously, $|\cP| \leq   \binom{n}{2} / \binom{r}{2}$.
If here equality holds, then $\cP$ is called an $S(n,r,2)$  {\em Steiner system}.

\begin{definition} A set system $ {\cal F}$ contains an $ a\times b$  {\em grid}  if
 there exist two disjoint subfamilies  ${\cal A}, {\cal B}\subseteq \cal F$ such that
\begin{itemize}
\item $|{\cal A}|=a$, $|{\cal B}|=b$, $|\cA\cup \cB|=a+b$,
\item $A\cap A'=B\cap B'=\emptyset$ for all $A,A'\in \cal A$, $A\neq A'$, $B,B'\in \cal B$, $B\neq B'$, and
\item $|A\cap B|=1$ for all $A\in \cA$, $B\in \cB$.
\end{itemize}
\end{definition}
Thus  an $r$-uniform $r\times r$  grid, $\Grr$, is a disjoint pair ${\cal A}, {\cal B}$
 of the same sizes $r$ such that they cover exactly the same set of $r^2$ elements.

\begin{theorem}\label{th:1}
For $r\ge 4$ there exists a real $c_r>0 $ such that there are linear $r$-uniform
 hypergraphs ${\cal F}$ on $n$ vertices  containing no grids  and
\begin{equation*}
|{\cal F}| >  \frac{n(n-1)}{r(r-1)}-c_rn^{8/5}.
\end{equation*}
\end{theorem}
The proof is postponed to Section~\ref{ss:32}.

The {\it Tur\'an number} of the $r$-uniform hypergraph $\cH$, denoted by $\ex(n,\cH)$,
 is the size of the largest $\cH$-free $r$-graph on $n$ vertices.
If we want to emphasize $r$, then we write $\ex_r(n, \cH)$.
Let $\I2$ be (more precisely $\I2^r$)  the class of hypergraphs of two edges
 and intersection sizes at least two.
This class consists of $r-2$  non-isomorphic hypergraphs, $\cI_j$, $2\leq j<r$,
 $\cI_j:=\{ A_j, B_j\}$ such that $|A_j|=|B_j|=r$, $|A_j\cap B_j|=j$.
Using these notations the above Theorem can be restated as follows.
 \begin{equation}\label{eq:2}
  \frac{n(n-1)}{r(r-1)}-c_rn^{8/5}<  \ex_r(n, \{ \I2, \Grr\}) \leq\frac{n(n-1)}{r(r-1)}
\end{equation}
holds for every $n,r\geq 4$.
In the case of $r=3$ we only have  
\begin{equation}\label{eq:r=3}
\Omega (n^{1.8})\leq   \ex_3(n, \{ \I2, \G33\}) \leq\frac{1}{6}n(n-1), 
  \end{equation}
  see in Section~\ref{ss:prob_lower}. 
The case of graphs, $r=2$, is different, see later in Section~\ref{ss:35}.

\begin{conjecture}\label{c:1}
The asymptotic {\rm (\ref{eq:2})} holds for $r=3$, too.\\ ---
Even more, for any given $r\geq 3$ there are infinitely many Steiner systems
  avoiding $\Grr$. \\ ---
Probably there exists an $n(r)$ such that, for every admissible $n>n(r)$
 (this means that $(n-1)/(r-1)$ and $\binom{n}{2}/\binom{r}{2}$ are both integers)
 there exists a grid-free $S(n,r,2)$.
  \end{conjecture}

\subsection{Sparse Steiner systems}

There are many problems and results concerning subfamilies of block designs, see, e.g., Colbourn and Rosa~\cite{CR}.
A Steiner triple system  ${\rm STS}(n):=S(n,3,2)$  is called $e$-{\em sparse} if
 it contains no set of $e$ distinct triples spanning at most $e+2$ points.
Every Steiner triple system is 3-sparse.
A longstanding conjecture of Erd\H os~\cite{ErdRoma} is that for every $e\ge 4$
 there exists an $n_0(e)$ such that if $n>n_0(e)$ and $n$ is admissible (i.e.,
  $n\equiv 1$ or $3$ (mod $6$)), then there exists an $e$-sparse ${\rm STS}(n)$.
Systems that are 4-sparse are those without a Pasch configuration (4 blocks spanning
 6 points, $\P4$, $\{ a,b,c\}$, $\{ a,d,e\}$, $\{b,d,f\}$, $\{c,e,f\}$).
Completing the works of Brouwer~\cite{Br},  Murphy~\cite{GM, GMP}, Ling and
 Colbourn~\cite{LCGG} and others, finally Grannell, Griggs and Whitehead~\cite{GGW}
 proved that 4-sparse ${\rm STS}(n)$'s
 exist for all admissible $n$ except 7 and 13.

   A 5-sparse system is precisely one lacking Pasch, $\P4$, and mitre configurations,
$\M5$, the latter comprising five blocks of the form
 $\{a,b,c\}$, $\{a,d,e\}$, $\{a,f,g\}$, $\{b,d,f\},$ $\{c,e,g\}$.
In a sequence of papers (e.g.,  Colbourn, Mendelsohn, Rosa, and
  \v Sir\'a\v n~\cite{CMRS}) culminating in  Y. Fujiwara~\cite{Fuy} and
  Wolfe~\cite{W06} it was established that systems having no mitres exist for
  all admissible orders, except for $n=9$.

Concerning the even more difficult problem of constructing 5-sparse systems (see Ling~\cite{Ling})
 Wolfe~\cite{W05, W08} proved that such systems exist for almost all admissible $n$.
More precisely, let  $A(x):=\{n\colon n\equiv 1$ or 3 (mod 6), $n\le x\}$ and
 $S(x):=\{n\colon$ there exists a 5-sparse ${\rm STS}(n)$ with $n\le x\}$, then
 $\lim_{x\to \infty}(|S(x)|/|A(x)|)=1$.

Forbes, Grannell and Griggs~\cite{FGG07, FGG09} constructed infinite classes of
 6-sparse ${\rm STS}(n)$'s.
As Teirlinck~\cite{T} writes  in his 2009 review of~\cite{W08}
``currently no nontrivial example of a 7-sparse Steiner triple system is known''.

Our Conjecture~\ref{c:1} is related to but not a consequence of Erd\H os' problem.
Colbourn~\cite{C} has checked (using a computer) all the 80 different
 ${\rm STS}(15)$'s and each contained at least 11 copies of ${\mathbb G}_{3\times 3}$.
Blokhuis~\cite{Bl} reformulated (a weaker version of) Conjecture~\ref{c:1} as follows:
Are there latin squares without the following subconfiguration?
$$
\left(
\begin{array}{ccc}
    *&a&b \\
     a&*& c\\
     b& c&*
   \end{array}\right)
   $$

Although the evidence is scarce one is tempted to generalize.
An $S(n,r,2)$ is $e$-{\em sparse} if the union of any $e$ blocks exceeds $e(r-2)+2$.

\begin{conjecture}\label{c:2}
 For every $e\ge 4$ there exists an $n_0(e,r)$ such that if
$n>n_0(e,r)$ and $n$ is admissible, then there exists an $e$-sparse $S(n,r,2)$.
  \end{conjecture}

\subsection{Sparse hypergraphs}

Brown, Erd\H os and  S\'os~\cite{E64, BES3, BESr} introduced the function $f_r(n,v,e)$ to denote
the maximum number of edges in an $r$-uniform hypergraph on $n$ vertices which does not contain
 $e$ edges spanned by $v$ vertices.
Such hypergraphs are called  $\bbG(v,e)$-{\em free} (more precisely $\bbG_r(v,e)$-{\em free}).
They showed that $f_r(n,e(r-k)+k,e)=\Theta(n^k)$ for every $2 \leq k < r$ and $e \geq 2$,
especially $f_r(n,e(r-2)+2,e)=\Theta(n^2)$.
The upper bound $\binom{n}{2}/\binom{r}{2}$ is easy, and this was the source of Erd\H os'
 conjecture concerning sparse Steiner systems. On the other hand, if
 we forbid $e$ edges spanned by one more vertex this problems
 becomes much more difficult.
Brown, Erd\H os and  S\'os conjectured that
\begin{equation}\label{eq:Econj}
   f_r(n,e(r-k)+k+1,e)=o(n^k).
 \end{equation}
One of the most famous results of this type is the $(6,3)$-Theorem of Ruzsa and Szemer\'edi~\cite{RSz},
the case $(e,k,r)=(3,2,3)$, saying that if no six points contain three triples then the size of the triple
system is $o(n^2)$, on the other hand  $n^{2-o(1)} < f_3(n,6,3)$.
This was extended by  Erd\H os, Frankl, and R\"odl~\cite{EFR} for arbitrary fixed $r\geq 3$,
\begin{equation}\label{eq:EFR}
n^{2-o(1)}<f_r(n,3(r-2)+3,3)=o(n^2).
  \end{equation}
The case $e=3$ was further extended by Alon and Shapira~\cite{AS}
 $$n^{k-o(1)} < f_r(n,3(r-k)+k+1,3) = o(n^k). $$
Even the case $k=2$,  $f_r(n,e(r-2)+3,e)=o(n^2)$,  is still open.
Nearly tight upper bounds were established by  S\'ark\"ozy and
Selkow~\cite{SS2, SSk}:
$$
  f_r(n,e(r-k)+k+\lfloor\log_2 e\rfloor,e)=o(n^k)\quad \forall r>k\geq2\ {\rm and}\ e\geq 3,
  $$
and for the case $e=4$, $r>k\geq 3$
$$
  f_r(n,4(r-k)+k+1,4)=o(n^k).
  $$

\subsection{Neither grids nor triangles}

\begin{definition}  Three sets $C_1,C_2,C_3$ form a  {\em triangle}, $\bbT_3$,
 if they pairwise intersect in three distinct singletons,
 $|C_1\cap C_2|=|C_2\cap C_3|=|C_3\cap C_1|=1$, $C_1\cap C_2\neq C_1\cap C_3$.
An $r$-uniform triangle is frequently denoted by $\bbT_3^r$.
\end{definition}

A {\em perfect matching} is a  subfamily  $\cM$ of the set system $ \cF$ such
 that the members of $\cM$ cover every element of $V(\cF)$ exactly once.

The main result of this paper is a construction.
\begin{theorem}\label{th:main}
For $r\ge 4$ there exist an $n_0(r)$ and $\beta_r>0$ such that
\begin{equation}\label{eq:th}
    \ex(n,\{ \I2, \bbT_3, \Grr \})  > n^2e^{-\beta _r\sqrt{\log n}}
  \end{equation}
holds for $n\ge n_0(r)$.
In other words, there exists a linear $r$-uniform hypergraph ${\cal F}$
 which contains neither grids nor triangles and $|{\cal F}|\ge n^2\, {\rm exp\, }[-\beta _r\sqrt{\log n}]$.
In addition, if $r$ divides $n$, then ${\cal F}$ can be  decomposed
  into perfect matchings, especially it is regular.

For the case $r=3$ we have the same with a much weaker lower bound
  \begin{equation}\label{eq:th_r=3}
    \ex(n,\{ \I2, \bbT_3, \G33 \})  > n^{1.6}e^{-\beta _3\sqrt{\log n}}. 
  \end{equation}

\end{theorem}
Again, the proof is postponed, to Section~\ref{ss:33}, and the cases $r\leq 3$
 to Section~\ref{ss:35}.

Note  that  $|{\cal F}|=o(n^2)$ by (\ref{eq:EFR}) so the lower bound (\ref{eq:th}) is almost optimal.
This result slightly improves the Erd\H os-Frankl-R\"odl~(\ref{eq:EFR})
 construction in two ways.
We make the hypergraph regular, and avoid not only triangles but grids, too.

\subsection{A probabilistic lower bound}\label{ss:prob_lower}

Almost all of the problems discussed in this paper can be formulated as a forbidden 
 substructure question, i.e., as  a Tur\'an type problem.
Here we present the standard probabilistic lower bound for the Tur\'an number
 due to Erd\H os, in a slightly stronger form as usual.
An $r$-uniform hypergraph $(V, \cF)$ is called $r$-{\em partite} if 
 there exists an $r$-partition of $V$, $V=V_1\cup \dots \cup V_r$, such that 
 $|F\cap V_i|=1$ for all $F\in \cF$, $i\in [r]$.
 
\begin{lemma}\label{le:erdos} {\rm (Erd\H os' lower bound on the Tur\'an number)}\newline
Suppose that $\cH$ is a (finite) family of $r$-graphs each of them having at least 
 two edges, and let
\begin{equation*} 
h:=\min\left\{ \frac{re-v}{e-1}: \bH\in \cH \text { is $r$-partite with $e$ edges 
 and $v$ vertices}\right\}.
  \end{equation*}
Then there exists a $c:=c(\cH)>0$ such that one can find an $n$-vertex $r$-partite
 $r$-graph $(n\geq r)$ of size at least $c n^h$ avoiding each member of $\cH$.
Hence
 \begin{equation}\label{eq:erdos_formula2}
    \ex(n,\cH)\geq \Omega (n^h). 
 \end{equation}
    \end{lemma}

\noindent
Sketch of the proof:
Choose independently each of the $(n/r)^r$ edges of the complete $r$-partite 
 hypergraph on $n$ vertices with probability $p$.
Leave out an edge from this random selection of each copy of $\bH\in \cH$.
The expected size of the remaining edges is at least 
$$
 p(n/r)^r -\sum p^e n^v. 
 $$
The rest is an easy calculation.  \qed

The ratio $(re-v)/(e-1)$ for $\I2^r$, $\Grr$, $\bbT_3$
 are $2$, $r^2/(2r-1)$, and $3/2$, resp., so Lemma~\ref{le:erdos} 
 implies (\ref{eq:r=3}), i.e., 
$\ex_3(\I2, \G33)\geq \Omega(n^{9/5})$. 
However, if $\bbT_3$ is among the forbidden substructures then 
 the probabilistic lower bound fails miserably, it gives only $\Omega(n^{3/2})$
 which is very far from the truth. 
To prove the slightly better lower bound (\ref{eq:th_r=3}) we are going to 
 use a version of the original Ruzsa-Szemer\'edi method.

\subsection{Tallying up the Tur\'an type problems}\label{ss:tally}

The three forbidden configurations,  $\I2^r$, $\Grr$, $\bbT_3^r$, have 7 non-empty 
  combinations. 
The cases $\ex(n,\{ \I2, \Grr \})$ and 
  $\ex(n,\{ \I2, \bbT_3, \Grr \})$ were discussed above in 
 Theorems~\ref{th:1} and~\ref{th:main}, 
 respectively. 
It is easy to see that 
\begin{equation*} 
\ex_r(n,\{ \I2^r, \bbT_3 \})= f_r(n,3(r-2)+3,3)+O(n),
   \end{equation*}
 so the Ruzsa and Szemer\'edi~\cite{RSz} and the Erd\H os, Frankl, and R\"odl~\cite{EFR} 
theorems, see (\ref{eq:EFR}), determine the right order of magnitude, $O(n^{2-o(1)})$. 

It was conjectured by Chv\'atal and Erd\H os and proved by Frankl et al.~\cite{Fur60} that 
\begin{equation*} 
\ex_r(n, \bbT_3 )= \binom{n-1}{r-1}
   \end{equation*}
for $r\geq 3$ and $n> n_0(r)$.
The only extremal $r$-graph consists of all $r$-tuples sharing a common element.
This hypergraph has no grid either, so we have
\begin{equation*} 
\ex(n,\{\bbT_3, \Grr \})= \binom{n-1}{r-1}
   \end{equation*}
for the same range of $r$ and $n$. 

We have $\ex_r(n, \I2^r)=\binom{n}{2}/\binom{r}{2}$ if and only if 
 a Steiner system $S(n,r,2)$ exists, which problem was solved for $n> n_0(r)$ 
 by Wilson~\cite{Wilson} and the exact packing number was determined for all 
 large enough $n$ by Caro and Yuster~\cite{CY}. 

The grid cannot be covered by $r-1$ vertices, it has $r$ disjoint edges.
So the $r$-graph having all edges meeting an $(r-1)$-element set is 
 grid free. 
This gives the  lower bound for the last case out of the seven.
\begin{equation}\label{eq:onlyGrr}
\ex(n, \Grr )\geq  \binom{n-1}{r-1}+ \binom{n-2}{r-1}+\dots +\binom{n-r+1}{r-1}. 
   \end{equation}
The classical result concerning the Tur\'an number of the complete $r$-partite 
 graph on $r\times r$ vertices by  Erd\H os~\cite{E64Isr} gives only an upper bound 
 $O(n^{r-\delta})$ with $\delta = r^{-r+1}$.
The truth should be much closer to the lower bound.

\begin{problem}
Determine the order of magnitude of $\ex(n, \Grr)$. 
  \end{problem}

\subsection{Union-free and cover-free hypergraphs}

Union free families were introduced by Kautz and Singleton
\cite{KS}. They studied binary codes with the property that the
disjunctions (bitwise $OR$s) of distinct at most $r$-tuples
of codewords are all different. In information theory
usually these codes are called {\em superimposed} and they have been
 investigated in several papers on multiple access communication (see,
 e.g., Nguyen Quang A and Zeisel~\cite{AZ}, D'yachkov and Rykov~\cite{DR1},
 Johnson~\cite{J1,J2,J3}).
Alon and Asodi~\cite{AA1,AA2}, and De Bonis
and Vaccaro~\cite{DV} studied this problem in a more general setup.
Small values of generalized superimposed codes and and their relation to
 designs were considered by  Kim, Lebedevin and Oh ~\cite{KL,KLO}.

The same problem has been posed -- in different terms --
by Erd\H os, Frankl and F\"uredi~\cite{EFF1, EFF2} in
combinatorics, by S\'os \cite{S} in combinatorial number theory, and
by Hwang and S\'os~\cite{H, HS} in group testing. One can
find short proofs of the best known upper bounds of these codes in the
papers by the present authors in~\cite{F} and~\cite{Ru}.
In \cite{FR} the connection of these codes to the big distance ones is shown.
A geometric version has been
posed by Ericson and Gy\"orfi \cite{EGy} and later investigated in~\cite{FuR}.
For a direct geometry application,
 notice that a union-free family defines a set of points of exponential size
in ${\bf R}^n$ such that arbitrary three of them span an acute triangle~\cite{EF}.

A family ${\mathcal F}\subseteq 2^{[n]}$ is $e$-{\em union-free}
if for arbitrary two distinct subsets ${\mathcal A}$ and ${\mathcal B}$ of
${\mathcal F}$ with $0< |{\mathcal A}|,|{\mathcal B}|\le e$
$$
 \bigcup_{A\in {\mathcal A}} A\not =\bigcup_{A\in {\mathcal B}} B.
$$
Let $U(n,e)$ ($U_r(n,e)$) be the maximum size of an $e$-union-free
 $n$ vertex hypergraph ($r$-uniform hypergraph, resp.).
The order of magnitude of $U_r(n,2)$ was determined by Frankl et al.~\cite{Fur25, Fur51}.

A family ${\mathcal F}\subseteq 2^{[n]}$ is $e$-{\em cover-free}
if for arbitrary distinct members $A_0,A_1,\dots , A_e\in {\mathcal F}$
$$A_0\not\subseteq\bigcup_{i=1}^e A_i.$$
Let $C(n,e)$ ($C_r(n,e)$) be the maximum size of an $e$-cover-free
 $n$ vertex hypergraph ($r$-uniform hypergraph, resp.).

An $e$-{cover-free} hypergraph is $e$-{union-free} and
an  $e$-{union-free} is $(e-1)$-{cover-free}.
(Indeed, the existence of an $(e-1)$-cover
$A_0\not\subseteq A_1\cup \dots \cup  A_{e-1}$ gives
 $\bigcup_{0\leq i\leq e-1}A_i=\bigcup_{1\leq i\leq e-1} A_i$).
Therefore,
\begin{equation}\label{eq:un}
    C(n,e)\leq U(n,e) \leq C(n,e-1)\leq U(n,e-1)\dots
  \end{equation}
and
\begin{equation}\label{eq:ur}
    C_r(n,e)\leq U_r(n,e) \leq C_r(n,e-1)\leq U_r(n,e-1)\dots
  \end{equation}
We have $C_r(n,r)=n-r+1$ (for $n\geq r$). Indeed,
 every member of an $r$-uniform $r$-cover-free family has a vertex of degree one.
In this section, based on Theorem~\ref{th:main},
 we determine the next two terms of the sequence (\ref{eq:ur}).
First, observe that
\begin{equation}\label{eq:8}
  C_r(n,r-1)\leq \frac{n(n-1)}{r(r-1)}.
  \end{equation}
Indeed, an $r$-uniform, $(r-1)$-cover-free family either has a vertex of degree one
 (and then we use induction on $n$), or it is  a linear hypergraph.

\begin{corollary}\label{co:15}
There exists a $\beta=\beta(r)> 0$ such that for all $n\geq r\geq 4$
 $$  n^2e^{-\beta _r\sqrt{\log n} }< U_r(n,r) \leq \frac{n(n-1)}{r(r-1)}.$$
In addition, if $r$ divides $n$, then our $n$-vertex, $r$-uniform,
 $r$-union-free family yielding the lower bound can be  decomposed
 into perfect matchings, especially it is regular.
  \end{corollary}
{\bf Proof.}
The upper bound follows from (\ref{eq:ur}) and (\ref{eq:8}), i.e.,
$$
  U_r(n,r) \leq C_r(n,r-1)\leq \frac{n(n-1)}{r(r-1)}.
  $$
On the other hand, we claim that
\begin{equation}\label{eq:conj9}
  \ex(n,\{ \I2^r, \bbT_3, \Grr \}) \leq U_r(n,r),
  \end{equation}
  hence the lower bound for $U_r(n,r)$ follows from (\ref{eq:th}).

We have to show that a linear $r$-graph without triangle and grid
 is $r$-union-free.
Suppose, on the contrary, that ${\cal A}\neq{\cal B}$,  $|{\cal B}|\le |{\cal A}|\le r $,
 $\cup_{A\in {\cal A}}A=\cup_{B\in {\cal B}}B$ and $\cA\cup \cB$
 form a linear $r$-uniform hypergraph.
Then $\exists A_1\in{\cal A}\setminus {\cal B}$.
Since $|A_1\cap B|\leq 1$, to cover the elements of $A_1$ the family ${\cal B}$
 must contain $r$ sets, i.e., $|{\cal B}|=|{\cal A}|= r$.
Moreover, the sets $B_1, ... , B_r\in \cB$ meet $A_1$ in distinct elements.
If ${\cB}$ consists of disjoint sets only, then $|\cup_{B\in {\cal B}} B|=r^2$, and to
cover these $r^2$ elements ${\cal A}$ must consist of disjoint
sets, too, and $\cA \cup \cB$ form a grid $\Grr$.
Otherwise, $\exists B_i, B_j\in {\cal B}$, such that $B_i\cap B_j=\{x\}\notin A_1$.
Then $A_1$, $B_i$, and $B_j$ form a triangle.
\qed

\begin{proposition}\label{pr:10}
In the case $r=3$ the probabilistic lower bound {\rm (\ref{eq:erdos_formula2})} implies 
\begin{equation}\label{eq:u3}
\Omega(n^{5/3})\leq U_3(n,3) 
  \end{equation}
\end{proposition}
The details are postponed to Section~\ref{ss:34}.

Let $\bbP_r$ be an $r$ uniform hypergraph with edges $A$, $B$ and 
 $C_1, \dots, C_{r-1}$ as follows.
The $r$-sets $C_1,\dots, C_{r-1}$ are pairwise disjoint, $a_i, b_i\in C_i$ are distinct elements,
 $d\notin \cup C_i$ and $A:=\{ d, a_1, a_2, \dots, a_{r-1}\}$ and $B:=\{ d, b_1, \dots, b_{r-1}\}$.

\begin{conjecture}\label{c:6}
If $\cF$ is an $n$-vertex, $r$-uniform ($r\geq 3$), linear hypergraph not containing $\bbP_r$, then
 its size $|\cF|=o(n^2)$. In other words,
$$
  \ex_r(n, \{ \I2, \bbP_r\})=o(n^2).
  $$
  \end{conjecture}
This would imply the conjecture of Erd\H os~(\ref{eq:Econj}) in the case
  $k=2$, $e=r+1$.
If it is true, then it implies the following more modest conjecture
\begin{equation}\label{eq:10}
 U_r(n,r)=o(n^2).
  \end{equation}

\begin{proposition}\label{pr:1}
Suppose that $r\ge 2$, $n\equiv r \pmod{ r^2-r}$, $1\leq k \leq (n-1)/(r-1)$
 and $n> n_0(r)$.
Then there exists a $k$-regular, $(r-1)$-cover-free $r$-graph.
Thus, in this case {\rm (\ref{eq:8})} gives
$$  C_r(n,r-1)=\frac{n(n-1)}{r(r-1)}.
  $$
\end{proposition}
{\bf Proof.}
To obtain the $k$-regular construction one can apply a classical theorem of
 Ray-Chaudhuri and Wilson~\cite{RWfelold}: For any given $r\geq 2$ there exists an
$n_0(r)$ such that, if $n>n_0(r)$ and $n\equiv r \pmod{r^2-r}$, then
 there exists a {\em resolvable}, $r$-uniform, $n$-vertex Steiner system $\mathcal S$.
This means that $\mathcal S$ can be decomposed into
$K:=(n-1)/(r-1)$ perfect matchings,
 (also called {\em parallel classes})
 ${\mathcal S}=\cup_{1\leq i\leq K}{\mathcal S}_i$, where
 $|{\mathcal S}_i|=n/r$ and $|\cup {\mathcal S_i}|=n$.
Taking $k$ of these parallel classes gives the desired
$(r-1)$-cover-free family. \qed

\subsection{Optimal superimposed codes}

D'yachkov and Rykov \cite{DR2} introduced the concept of
 {\em optimal superimposed codes and designs}.
Recall an easy observation.

\begin{proposition} {\rm (D'yachkov, Rykov \cite{DR2})}
\label{opt}
\quad
If ${\mathcal F}\subseteq 2^{[n]}$ is $(r-1)$-{cover-free}, $(r\geq 2)$, it has maximum degree $k$ and
$|{\mathcal F}|=t\ge n$ then $\left\lceil{nk/r}\right\rceil\geq t$ holds.
\end{proposition}
By (\ref{eq:un}) a similar statement holds for $r$-{union-free} families, too.
Note that, from coding theory point of view, it is reasonable to assume that $t\ge n$
since a collection of singletons ${\mathcal F}$ is $(r-1)$-{cover-free} for arbitrary
$2\le r\le n$ with $|{\mathcal F}|=n$ and usually the goal is to get a code as large
 as possible.

{\bf Proof of \ref{opt}.}\quad
Let ${\mathcal F}_0=\{A\in{\mathcal F}: \exists x\in A, \deg_\cF(x)=1\}$,
 $|{\mathcal F}_0|=t_0$.
Clearly, $|A|\ge r$ for every $A\in {\mathcal F}\setminus
  {\mathcal F}_0$ otherwise the union of some other $(r-1)$
  members of ${\mathcal F}$ cover $A$.
We obtain
\begin{equation*}
   r(t-t_0)+t_0\le \sum_{A\in \cF} |A| =\sum_{x\in [n]} \deg (x)\leq  k(n-t_0)+t_0 \, . \qed
 \end{equation*}

It follows that in case of $nk=rt$ the family $\cF$ should be $k$-regular and $r$-uniform.

\begin{definition}\label{def:opti} {\rm  (see~\cite{DR2})}\quad
The  $n$-vertex family $\cF$ is called an  {\em optimal $(r-1)$-superimposed code}
  if it is an $(r-1)$-cover-free, $r$-uniform, $k$-regular, linear hypergraph.
It is called an {\em optimal $r$-superimposed design} if in addition it is
 $r$-{union-free}, too.
In both cases $nk=rt$ holds.
\end{definition}

Let $k(r-1,n)$  $(k'(r,n))$ denote the maximum $k$ that such a $k$-regular optimal
 $(r-1)$-superimposed code (optimal $r$-superimposed design) exits.
They have showed for every $r\geq 2$ and $n$ that
\begin{equation*}
\begin{array}{rclccccl}
  {} & {}& k'(r,n)&{}&\leq& k(r-1,n) &\leq& (n-1)/(r-1)\\
\log_2 n -O(1)&\leq& k'(2,n), \quad\quad&  n/2 &\leq&  k(1,n)&\leq& \enskip n-1\\
  4&\leq& k'(3,n),\quad\quad &(n/3)-1&\leq &k(2,n)&\leq &(n-1)/2.
  \end{array}
  \end{equation*}
They and Macula~\cite{Macu} gave a lower bound for every $r\ge 3$ for a few special
 but infinitely many  values of $n$.
  \begin{equation}\label{DyR}
 \left(\frac{n}{r}\right)^{1/(r-1)}\leq k'(r,n).
 \end{equation}
Our results, Corollary~\ref{co:15} and Proposition~\ref{pr:1},  imply that

\begin{corollary}\label{co:th32:fo}
If $r\ge 4$, $r|n$ and $n\ge n_0(r)$, then there exists an optimal $k$-regular,
$r$-superimposed design for every $1\leq k\leq  n e^{-\beta _r\sqrt{\log n}}$,
 especially\quad
  $$ n e^{-\beta _r\sqrt{\log n}}\leq k'(r,n). \qed $$
\end{corollary}

\begin{proposition}\label{pr:16}
For the case $r=3$ we have the same statement with a weaker lower bound
\begin{equation}\label{eq:k3}
  \frac{1}{25}n^{2/3}\leq k'(3,n).  
  \end{equation}
\end{proposition}
The details are postponed to Section~\ref{ss:34}.

\begin{corollary}\label{co:th33}
Suppose that $r\ge 2$, $n\equiv r \pmod{ r^2-r}$ and $n> n_0(r)$.
Then there exists an  optimal  $(r-1)$-superimposed code for every
 $1\leq k\leq (n-1)/(r-1)$, especially $$k(r-1,n)=(n-1)/(r-1). \qed $$
\end{corollary}

\section{Tools from combinatorial number theory and discrete geometry}

\subsection{Three lemmata from combinatorial number theory}

\begin{lemma}\label{le:mink} {\rm (Minkowksi's theorem of simultaneous approximation \cite{M})}\quad
Let $q$ be a prime and $(n_1, \dots, n_d)\in R^d$ an integer point.
Then there exist an integer $0< \alpha < q$ and  residues $r_i$
 such that $r_i\equiv \alpha n_i\,  (\!\!\! \mod q)$ and  $|r_i|\leq q^{1-1/d}$
 for all $1\leq i \leq d$.
  \end{lemma}

\noindent
Sketch of the proof:
Consider all vectors of the form $a {\mathbf n}$ mod $q$, $a=0, 1,2, \dots, q-1$.
There will be two of them $a_1{\mathbf n}$ and $a_2{\mathbf n}$ `close' to each other.
Take $\alpha =a_1-a_2$. \qed

Let $r_k(q)$ be the  maximum number of
integers which can be selected from $\{1,\dots ,q\}$ containing no
$k$-term arithmetic progression.
This function has been extensively studied in the last six decades by leading
 mathematicians see, e.g., Ruzsa~\cite{ruzsa}.
The major important known bounds (appart from some recent minor improvements)
 are due to Behrend~\cite{Be}, Heath-Brown~\cite{HB} and Szemer\'edi~\cite{Sz}:
 there are positive constants $\alpha$ and $\beta$ such that
\begin{equation}\label{eq:Sz}
     qe^{-\beta\sqrt{\log q}}<r_3(q)<q(\log q)^{-\alpha}\quad\quad
     \text{\rm and for all $k$}\quad  r_k(n)=o(n).
     \end{equation}
Call a set $M\subset [q]$  $r$-{\em sum-free} if the equation
\begin{equation*} 
     c_1m_1 +c_2m_2=(c_1+c_2)m_3
     \end{equation*}
has no solutions with $m_1, m_2, m_3\in M$ and $c_1, c_2$  are positive integers
 with $c_1+c_2\le r$  except the one with $m_1=m_2=m_3$.
We will need the following lower bound used by Erd\H os, Frankl
and R\"odl~\cite{EFR}, also see Ruzsa~\cite{ruzsa}.
Its proof requires only a  slight modification of Behrend's~\cite{Be} argument.

\begin{lemma}\label{le:behrend} {\rm (Generalized Behrend)}\quad
For arbitrary positive integer $r$ there exists a $\gamma _r>0$
 such that for any integer $q$ one can find an $r$-sum-free subset
  $M\subseteq \{0,1,..., q \}$ such that  $|M|>qe^{-\gamma_r\sqrt{\log q}}$.
  \end{lemma}

The case $r=2$ (and $c_1=c_2=1$) is the original statement of
Behrend~\cite{Be}.
Ruzsa also notes that an upper bound $O(q/ (\log q)^{\alpha_r})$ for the general case
 can be proved by the methods of~\cite{HB} and~\cite{Sz}.

Call a set of numbers  {\em $A_6$-free} if it does not contain a subset
 of the form 
$$
 \{ x-a-b, x-b, x-a, x+a, x+b, x+a+b\}$$
 for some $a,b>0$, $a\neq b$.
Call it $A_4${\em -free} if it does not contain a fourtuple of the form 
 $\{ x-2a, x-a, x+a, x+2a\}$ for some $a>0$,
 and call it $AP_k${\em -free} if it contains no
 $k$-term arithmetic progression.
Let $r(n, P_1, P_2, \dots )$ denote the  maximum number of
integers which can be selected from $\{1,\dots ,n\}$ avoiding the patterns
$P_1, P_2, \dots$.
With this notation $r_3(n):= r(n,AP_3)$. 

Since an $A_4$-free set has no 5-term arithmetic progression we get 
 $r(n, A_4)\leq r_5(n)=o(n)$ by Szemer\'edi's~\cite{Sz} theorem~(\ref{eq:Sz}).
A 4-sum-free sequence is $A_4$-free as well (one has, e.g., 
 $1\times (x-2a)+ 3\times (x+2a)=4\times(x+a)$), thus
 Lemma~\ref{le:behrend} gives a lower bound showing
\begin{equation*}  
       r(n, A_4)=n^{1-o(1)}.  
  \end{equation*}

\begin{lemma}\label{le:A6}
\quad
\begin{equation*}  
      \frac{2}{5}r_3(n)^{3/5} < r(n, A_6, A_4, AP_3).  
  \end{equation*}
  \end{lemma}

\noindent
Sketch of the proof:
Similar to the proof of Lemma~\ref{le:erdos}.
Take an $AP_3$-free set $M\subset \{ 1, 2, \dots, n\}$ of maximum size.
Choose independently each element of $M$ with probability $p$, 
 and leave out an element from this random selection of each copy of
 the arising configurations we want to avoid.
The expected size of the remaining elements is at least 
$$
 p|M| - p^6|M|^3-p^4|M|^2. 
 $$
Define $p$ as $\frac{1}{2}|M|^{-2/5}$.  \qed

Starting with $M=[n]$, the same process gives
\begin{equation*}
   \frac{2}{5} n^{3/5} < r(n, A_6, A_4)\leq r(n, A_6).   
\end{equation*}
Concerning the upper bounds we only have $r(n, A_6)\leq r_7(n)=o(n)$. 

The random method notoriously gives a weak lower bound of Sidon type problems 
 (for definitions, see, e.g., Babai and S\'os \cite{BS}), e.g., the above argument 
 gives only $r(n, {\rm Sidon})\geq \Omega(n^{1/3})$, although the truth is 
 $\Theta(n^{1/2})$ (Erd\H os and Tur\'an~\cite{ET}) 
 and it gives $r_3(n, AP_3)\geq \Omega(n^{1/2})$ although the truth is
 $n^{1-o(1)}$. 
So we believe that there are much better lower bounds.
\begin{conjecture}
There is an $\varepsilon>0$ such that
  \begin{equation*}  
  n^{3/5 +\varepsilon} <r(n,A_6) 
   \end{equation*}
holds for large enough $n$.
Possibly the order of magnitude of this function is $n^{1-o(1)}$. 
  \end{conjecture}

\subsection{Grids of pseudolines}

A set of {\em pseudolines} is usually a set of (infinite) planer curves pairwise meeting
in at most one point with crossing, no two pseudolines are tangent.
The main result of this subsection is in fact deals with pseudoline arrangements
but we formulate it in a simpler way.

Let $\ell_1, \dots, \ell_r$ be parallel vertical lines on the plane,
 $\ell_j:=\{ (x,y): x=j\}$, $r\geq 2$.
Let $V_j$ be an $r$-set of points on the line $\ell_j$,
 $V_j:=\{ Q_{1,j}, \dots, ,Q_{r,j}\}$, $Q_{i,j}=(j, y_{i,j})$, such that
 $y_{1,j}>y_{2,j}>...>y_{r,j}$ for every $1\le j\le r$.
The points of $\Pi:= \cup V_j$ can be arranged in a matrix form
\begin{equation*}
\left(
\begin{array}{cccc} 
Q_{1,1} & Q_{1,2} & \dots & Q_{1,r}\\
Q_{2,1} & Q_{2,2} & \dots & Q_{1,r}\\
\vdots & \vdots & \ddots & \vdots \\
Q_{r,1} & Q_{r,2} & \dots & Q_{r,r}\\
\end{array}
\right)
\end{equation*}
where the elements of the $j$th column lie on $\ell_j$ and are
ordered from the top to the bottom. The vertical distances
$y_{i,j}-y_{i+1,j}$ are positive, but can be distinct from each
other.

A $\Pi$-{\em polygon} $\pi$ consists of $r-1$ segments of the form
  $\pi=\cup_j [Q_jQ_{j+1}]$ ($1\leq j\leq r-1$) with $Q_j\in V_j$.
There are $r^r$ such polygonal arcs.
Such a $\pi$ can be considered as a piecewise linear continuous function
  $\pi: [1,r]\to\bR$.
Two sets of $\Pi$-polygons $\cP$ and $\cR$ are called {\em crossing} if
\newline\indent (C1)\quad
 $|\cP|=|\cR|=r$,
\newline\indent (C2)\quad
 $\cP$ is covering all vertices of $\Pi$, i.e.,
$\{ \pi\cap \ell_j: \pi \in \cP\} =V_j$
and the same holds for $\cR$,
\newline\indent (C3)\quad
$\cP\cup\cR$ is almost disjoint, i.e.,
 $\pi, \rho\in (\cP\cup \cR)$ (and $\pi\neq \rho$) imply that $|\pi\cap \rho|\leq 1$, and
\newline\indent (C4)\quad
$\cP\cup\cR$ behaves like pseudolines, i.e.,
$\pi, \rho\in (\cP\cup \cR)$, $|\pi\cap \rho|=1$ imply that
they are truly crossing, i.e., if $\pi\cap \rho=(x_0,y_0)$ then
 \newline\indent either
$\pi(x)< \rho (x)$ for all $1\leq x< x_0$ together with
$\pi(x)> \rho (x)$ for $x_0< x\leq r$,
\newline\indent
or $\pi(x)> \rho (x)$ for all $1\leq x< x_0$ together with
$\pi(x)< \rho (x)$ for $x_0< x\leq r$.

Note that (C1) and (C2) imply that every member of $\Pi$ belongs to a unique
  $\pi \in \cP$ and also to a unique $\rho\in \cR$.
Then (C3) yields that each $\pi \in \cP$ meets each $\rho \in \cR$ at an element
 of $\Pi$, and only there.
The members of $\cP$ (and $\cR$) can be disjoint among themselves.

A {\em very special} crossing system is depicted below in (\ref{pattern}) using
 thin and thick lines to indicate the segments of ${\cP}$ and ${\cR}$.
Its properties described in (VS1)--(VS7) below.

\begin{equation} \label{pattern}
\begin{array}{cccccccc}
Q_{1,1} &
\mbox{\begin{picture}(35,11)\thinlines\put(0,3){\line(1,0){35}}
\end{picture}}
& Q_{1,2} &
\mbox{\begin{picture}(35,11)\thicklines\linethickness{2pt}\put(0,3){\line(1,0){35}}
\end{picture}}
& Q_{1,3} &
\mbox{\begin{picture}(35,11)\thinlines\put(0,3){\line(1,0){35}}
\end{picture}} & Q_{1,4}
\\
~ &
\mbox{\begin{picture}(35,35)\thicklines\linethickness{2pt}\put(0,0){\line(6,5){35}}
\put(35,0){\line(-6,5){35}}\end{picture}} & ~ &
\mbox{\begin{picture}(35,35)\thinlines \put(0,0){\line(6,5){35}}
\put(35,0){\line(-6,5){35}}\end{picture}} & ~ &
\mbox{\begin{picture}(35,35)\thicklines\linethickness{2pt}
\put(0,0){\line(6,5){35}}
\put(35,0){\line(-6,5){35}}\end{picture}} & ~ \\
Q_{2,1} & ~ & Q_{2,2} & ~ & Q_{2,3} &~& Q_{2,4} &~ \\
~ & \mbox{\begin{picture}(35,35)\thinlines \put(0,0){\line(6,5){35}}
\put(35,0){\line(-6,5){35}}\end{picture}} & ~ &
\mbox{\begin{picture}(35,35)\thicklines\linethickness{2pt}
\put(0,0){\line(6,5){35}} \put(35,0){\line(-6,5){35}}\end{picture}}
& ~ & \mbox{\begin{picture}(35,35)\thinlines
\put(0,0){\line(6,5){35}}
\put(35,0){\line(-6,5){35}}\end{picture}}&~&\raisebox{2.8ex}{$\cdots$}
\\
Q_{3,1} & ~ & Q_{3,2} & ~ & Q_{3,3} &~& Q_{3,4} & ~ \\
~ & \mbox{\begin{picture}(35,35)\thicklines\linethickness{2pt}
\put(0,0){\line(6,5){35}} \put(35,0){\line(-6,5){35}}\end{picture}}
 & ~ &
\mbox{\begin{picture}(35,35)\thinlines \put(0,0){\line(6,5){35}}
\put(35,0){\line(-6,5){35}}\end{picture}}
 & ~ &
\mbox{\begin{picture}(35,35)\thicklines\linethickness{2pt}
\put(0,0){\line(6,5){35}}
\put(35,0){\line(-6,5){35}}\end{picture}}&~\\
Q_{4,1}&~&Q_{4,2}&~ &Q_{4,3}&~&Q_{4,4}\\
~&~&~&\vdots&~&~&~&\ddots
\end{array}
\end{equation}

There are two types of segments defining $\cP\cup \cR$,
\newline\noindent
(VS1)\quad
 all edges of the upper and lower envelops $[Q_{1,j}, Q_{1,j+1}]$,
 $[Q_{r,j}, Q_{r,j+1}]$, $1\leq j\leq r-1$, and \newline
(VS2)\quad all diagonal edges $Q_{i,j}Q_{i\pm1,j+1}$. \newline
(VS3)\quad The edges on the upper (lower) envelops alternate
 between $\cP$ and $\cR$. \newline
(VS4)\quad The crossing diagonal edges $[Q_{i,j},Q_{i+1,j+1}]$ and $[Q_{i+1,j},Q_{i,j+1}]$
 simultaneously belong to either ${\cP}$ or ${\cR}$. \newline
\indent $\cP\cup \cR$ consists of the following three types of polygonal paths:
\newline
(VS5)\quad
  the two diagonals $Q_{1,1}, Q_{2,2}, \dots, Q_{r,r}$, and $Q_{r,1}, Q_{r-1,2}, \dots, Q_{1,r}$,\newline
(VS6)\quad for $1\leq j \leq r-1$ a {\em top} path consisting of an increasing part (on Figure {\rm (\ref{pattern})})
  $Q_{j,1}$, $Q_{j-1,2}, \dots, ,Q_{1,j}$ a horizontal edge $Q_{1,j}, Q_{1,j+1}$ and a decreasing part
 $Q_{1,j+1}$, $Q_{2,j+2}, \dots, Q_{r-j,r}$,\newline
(VS7)\quad  for each $2\leq j \leq r$ a {\em bottom} path starting at $Q_{j,1}$ and
  consisting of a decreasing part having vertices $Q_{j+x,1+x}$ $(x=0,1,2,\dots, r-j)$,
 a horizontal edge $[Q_{r,r-j+1}, Q_{r,r-j+2}]$ and an
 increasing part
  $Q_{r-x,r-j+2+x}$ $(x=0,1,2,\dots, j-2)$.

\begin{lemma}\label{le:cross}
Suppose that the two sets of $\Pi$-polygons $\cP$
 and $\cR$ form a crossing $\Pi$-polygon system (i.e., satisfy {\rm (C1)--(C4)}).
Then they have the unique structure described by {\rm (VS1)--(VS7)}.
  \end{lemma}

The proof is postponed to the next subsection.
From this Lemma we can read out the intersection structure. We obtain

\begin{equation}\label{eq:8x8}
\includegraphics[width=8cm]{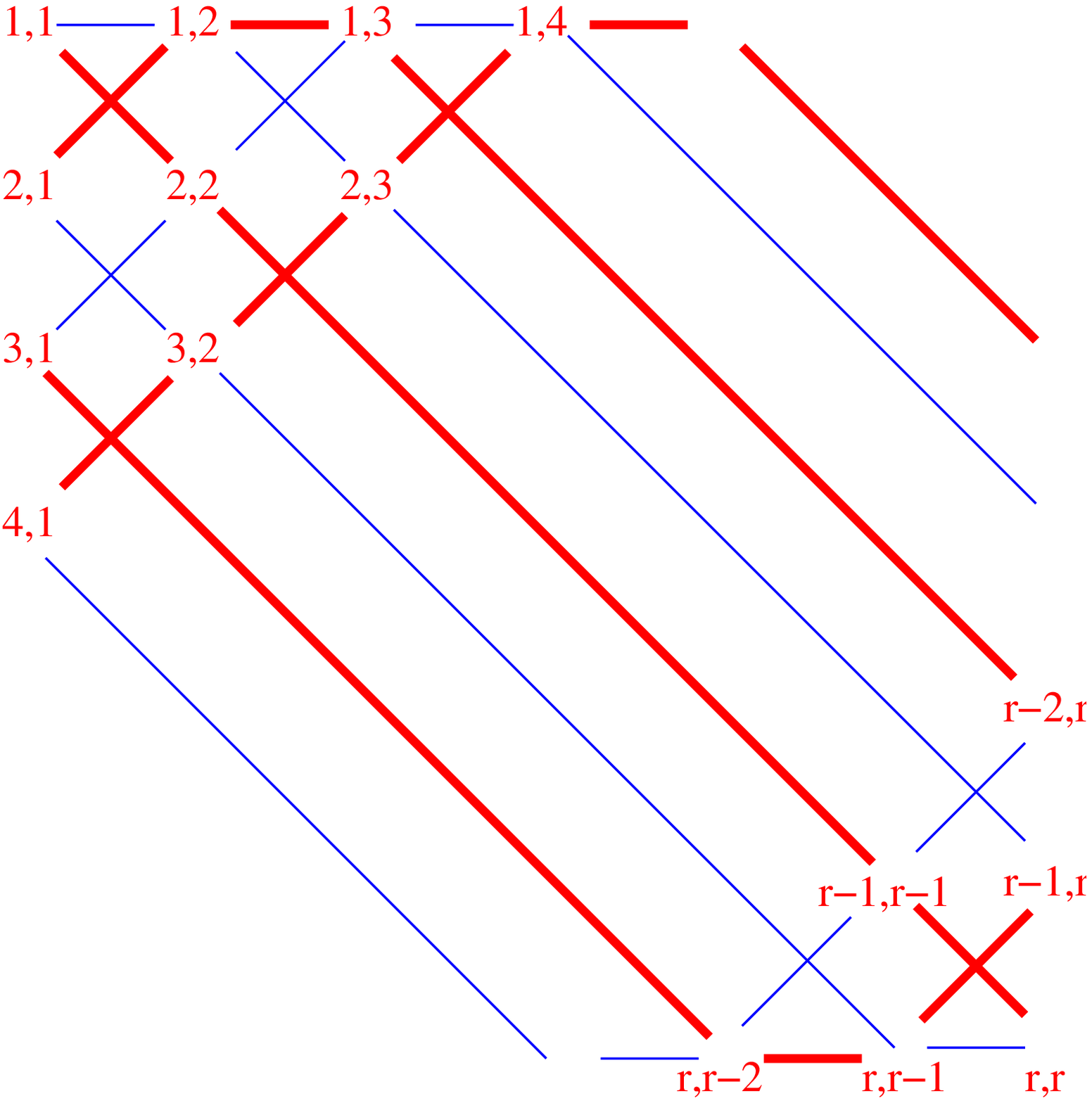}
\end{equation}

\begin{corollary}\label{co:1-4}
Suppose that the two sets of $\Pi$-polygons $\cP:=\{ \pi_1, \dots, \pi_r\}$
 and $\cR:=\{ \rho_1, \dots, \rho_r\}$ form a crossing $\Pi$-polygon system
 with $Q_{i,1}\in \pi_i, \rho_i$ $(1\leq i\leq r)$ and with $[Q_{1,1},Q_{1,2}]\subset \pi_1$.
Then
\newline --- the vertices of $\pi_1$ are $Q_{1,1}, Q_{1,2}, Q_{2,3}, \dots,   Q_{r-1,r}$,
\newline --- the vertices of $\pi_2$ are $Q_{2,1}, Q_{3,2},  \dots, Q_{r,r-1}, Q_{r,r}$,
\newline --- the vertices of $\pi_3$ are $Q_{3,1}, Q_{2,2}, Q_{1,3}, Q_{1,4}, \dots, Q_{r-3,r}$,
\newline --- the vertices of $\pi_4$ are $Q_{4,1},   \dots, Q_{r,r-3}, Q_{r,r-2}, Q_{r-1,r-1}, Q_{r-2,r}$,
\newline --- the vertices of $\rho_1$ are $Q_{1,1}, Q_{2,2},  \dots, Q_{r-1,r-1}, Q_{r,r}$,
\newline --- the vertices of $\rho_2$ are $Q_{2,1}, Q_{1,2}, Q_{1,3},  \dots,  Q_{r-2,r}$,
\newline --- the vertices of $\rho_3$ are $Q_{3,1}, \dots, Q_{r,r-2},Q_{r,r-1}, Q_{r-1,r}$,
\newline --- the vertices of $\rho_4$ are $Q_{4,1}, Q_{3,2}, Q_{2,3}, Q_{1,4}, Q_{1,5}, \dots, Q_{r-4,r}$. \qed
  \end{corollary}

\subsection{The proof of the uniqueness of the crossing structure}

Here we prove Lemma~\ref{le:cross} with a series of propositions.
Assume $Q_{1,i}=\pi_i\cap \rho_i$ for $1\leq i\leq r$.

For $1\le j \le r-1$ let $G_j^{\cP}$ be the bipartite graph (a matching) with
 parts $V_j$ and $V_{j+1}$ and with edges defined
 by the corresponding parts of the polygons from $\cP$, i.e.,
 $[Q_{i,j},Q_{k,j+1}]\in E(G_j^{\cP})$ if and only if 
there is a $\pi\in \cP$ with $\pi(j)=y_{i,j}$ and $\pi(j+1)=y_{k,j+1}$.
(Thus, to simplify notations, we identify the graph $G_j^{\cP}$ with its
 geometric representation.)
 $G_j^{\cR}$ is defined similarly. Finally,
 $G^{\cP}$ is the union of $G_j^{\cP}$ and the graph $G$
 is having all the edges of $G^{\cP}$ and $G^{\cR}$.

\begin{proposition}\label{egy}
Suppose $\pi\cap \ell_j=Q_{i,j}$ for some $\pi \in \cP$.
Then $\pi\cap\ell_{j+1}\in\{Q_{i-1,j+1},Q_{i,j+1},Q_{i+1,j+1}\}$.
Similarly, $\rho\in \cR$, $1\leq j\leq r-1$, $\rho\cap\ell_j=Q_{i,j}$ and
 $\rho\cap \ell_{j+1}=Q_{k,j+1}$ imply $|i-k|\leq 1$.
\end{proposition}
{\bf Proof.}\quad
We prove the second statement.
Assume to the contrary that $k\leq i-2$ (the case $k\geq i+2$ is similar).
Consider the $i-1$ edges of $G^\cP_j$ with vertices $Q_{1,j}, \dots, Q_{i-1,j}$.
Since there is not enough room to match these vertices to
$Q_{h,j+1}$ $(1\le h\le k)$  there exits a
 $[Q_{u,j},Q_{v,j+1}]\in E(G^\cP_j)$ with
 $u<i$ and $v>k$.
Then this segment intersects
$[Q_{i,j},Q_{k,j+1}]\in E(G_j^\cR)$ inside the open strip
 $\{ (x,y): j< x< j+1\}$. This contradicts to
the fact that a $\pi\in \cP$ and a $\rho\in \cR$ meet only in
 the points of $\Pi$. \qed

In the same way we obtain the following.

\begin{proposition} \label{pr:endpoint}
Suppose that $\gamma\in \cP\cup\cR$, the
 first point of $\gamma$ is $Q_{a,1}$, the last one is $Q_{b,r}$.
Then $b\in \{ r-a, r-a+1, r-a+2 \}$.
\end{proposition}
{\bf Proof.}\quad
Suppose $\gamma=\pi_a$.
Every $\rho_1,\dots,\rho_{a-1}$ starts above $\pi_a$ on $\ell_1$,
 they meet $\pi_a$ in a point of $\Pi$, so their endpoints
 on $\ell_r$ lie below or on the endpoint of $\gamma$, $Q_{b,r}$.
Hence $b\leq r-(a-2)$.
Considering $\rho_{a+1}, \dots, \rho_r$ the same argument gives $b\geq r-a$. \qed

Now we prove Lemma~\ref{le:cross} by checking (VS1)--(VS7).

\medskip\noindent{\bf Proof of (VS5).}\quad
The paths $\pi_1$ and $\rho_1$ end at either $Q_{r-1,r}$ or at $Q_{r,r}$ by the
 previous Proposition. They cannot meet in a second point, so
 one of them finishes at $Q_{r,r}$.
Then Proposition~\ref{egy} implies that this path is the diagonal of the
 $r\times r$ array,
 $Q_{1,1}, Q_{2,2}, Q_{3,3},  \dots, Q_{r-1,r-1}, Q_{r,r}$.
By symmetry, the other diagonal
  $Q_{r,1}, Q_{r-1,2},  \dots, Q_{1,r}$ also belongs to $\cP\cup \cR$.\qed

\medskip\noindent{\bf Proof of (VS1).}\quad
We show that each top edge, $[Q_{1,j},Q_{1,j+1}]$, belongs to $E(G)$.
According to Proposition~\ref{egy}
 the neighbor of $Q_{1,j}$ in $G_j^\cP$ (in $G_j^\cR$) is
 either $Q_{1,j+1}$ or $Q_{2,j+1}$.
The $\cP$ and $\cR$ edges are distinct, so both
 $[Q_{1,j},Q_{1,j+1}]$ and $[Q_{1,j},Q_{2,j+1}]$
must belong to $E(G)$.
Similar argument gives that both
 $[Q_{1,j},Q_{1,j+1}]$ and $[Q_{2,j},Q_{1,j+1}]$
 belong to $E(G)$ and we got the edges of the top layer of (\ref{pattern}). \qed

\medskip\noindent{\bf Proof of (VS3).}\quad
We show that the top edges alternate between $G^\cP$ and $G^\cR$.
If two consecutive of them, $[Q_{1,j},Q_{1,j+1}]$ and $[Q_{1,j+1},Q_{1,j+2}]$,
 both belong to $G^\cP$ then they are part of the same $\pi\in \cP$
and the diagonal edges $[Q_{1,j},Q_{2,j+1}]$ and $[Q_{2,j+1},Q_{1,j+2}]$
 are necessarily $G^\cR$ edges, belonging to the same  $\rho\in \cR$.
Then $\pi$ and $\rho$  cross twice, violating (C3). \qed

\medskip\noindent{\bf Proof of (VS6).}\quad
Call a path $\gamma\in \cP\cup\cR$ a {\em top} ({\em bottom}) path if it contains
 a top edge $[Q_{1,j}, Q_{1,j+1}]$ (bottom edge  $[Q_{r,j}, Q_{r,j+1}]$).
The pseudoline structure implies that these paths are all distinct, so $\cP\cup \cR$
 consists of the two diagonals and $r-1$ top and $r-1$ bottom paths.
Again Proposition~\ref{egy} implies that if $[Q_{1,j}, Q_{1,j+1}]\subset \gamma$ and the
 first point of $\gamma$ is $Q_{a,1}$, the last is $Q_{b,r}$ then $a\leq j$ and
 $b\leq r-j$ so $a+b\leq r$.
This and Proposition~\ref{pr:endpoint} give that $b=r-a$ and $\gamma$ must have the
 shape as described in (VS6).
By symmetry, we have the same for the bottom paths, property (VS7).
  \qed

Properties (VS6--7) imply that all diagonal edges $[Q_{i,j},Q_{i\pm1,j+1}]$ are in $E(G)$,
 this is property (VS2).
Property (VS4) follows from (VS2) and from the fact that a $\cP$ and an $\cR$ edge can meet
 only at the vertices of $\Pi$.
This completes the proof of Lemma~\ref{le:cross}.  \qed

\subsection{Grids of Euclidean lines}

Here we apply the results of the previous subsections when the
sets of crossing polygons $\cP$ and $\cR$ are actually straight
lines.
Let $\mathbf  y$, $\mathbf  m$, ${\mathbf  m}^\prime\in \bR^r$,
 ${\mathbf  y}=(y_1, \dots, y_m)$ with $y_1> y_2> \dots, y_r$,
 $\pi_i:=\{ (1,y_i), (2,y_i+m_i), \dots, (r,y_i+(r-1)m_i)\}$,
 $\rho_i:=\{ (j,y_i+(j-1)m_j^\prime): 1\leq j\leq r\}$,
 $\cP:=\{ \pi_1, \dots, \pi_r\}$, $\cR:=\{ \rho_1, \dots, \rho_r\}$
and $V_j:=\{\ell_j\cap \pi_i: i\in [r]\}= \{\ell_j\cap \rho_i: i\in [r] \}$
 with $\Pi:=\cup V_j$, $|\Pi|=r^2$.

\begin{lemma}\label{le:28}
If $\cP$ and $\cR$ are forming a crossing pair of straight lines, then $r\leq 3$.
  \end{lemma}

\noindent{\bf Proof.}\quad
$\cP$ and $\cR$ satisfy (C1)--(C4) so Lemma~\ref{le:cross} can be applied.
Consider the first four lines $\pi_1, \dots, \pi_4\in \cP$ and
  also $\rho_1,\dots, \rho_4\in \cR$.
According to Corollary~\ref{co:1-4} these lines meet
 at $Q_{1,1}$, $Q_{2,1}$, $Q_{3,1}$, and $Q_{4,1}$ and
 at 12 further points of $\Pi$ (see (\ref{eq:8x8})).
These 12 points yield 12 equations for $y_1, \dots, y_4, m_1,
\dots, m_4, m_1', \dots, m_4'$, e.g.,
 considering $Q_{1,2}$ we get the equation $y_1+m_1=y_2+m_2'$.
So the points $Q_{1,2}$, $Q_{1,3}$, $Q_{1,4}$, $Q_{2,2}$, $Q_{2,3}$, $Q_{3,2}$, and
  $Q_{r-2,r}$, $Q_{r-1,r-1}$, $Q_{r-1,r}$, $Q_{r,r-2}$, $Q_{r,r-1}$, and $Q_{r,r}$
   give the following 12 equations.
\begin{equation}\label{eq:12}
\left(
\begin{array}{cccccccccccc}
1&-1&0&0  \enskip &1&0&0&0  \enskip &0&-1&0&0\\
0&-1&1&0  \enskip &0&0&2&0   \enskip &0&-2&0&0\\
0&0&1&-1  \enskip &0&0&3&0  \enskip &0&0&0&-3\\
-1&0&1&0  \enskip &0&0&1&0  \enskip &-1&0&0&0\\
1&0&0&-1  \enskip &2&0&0&0  \enskip &0&0&0&-2\\
0&1&0&-1  \enskip &0&1&0&0  \enskip &0&0&0&-1
\\
\\
0&-1&0&1 \enskip  &0&0&0&r-1  \enskip &0&-r+1&0&0\\
-1&0&0&1  \enskip &0&0&0&r-2 \enskip &-r+2&0&0&0\\
1&0&-1&0  \enskip &r-1&0&0&0  \enskip &0&0&-r+1&0\\
0&0&-1&1  \enskip &0&0&0&r-3  \enskip &0&0&-r+3&0\\
0&1&-1&0  \enskip &0&r-2&0&0  \enskip &0&0&-r+2&0\\
-1&1&0&0  \enskip &0&r-1&0&0 \enskip  &-r+1&0&0&0
\end{array}\right)
\left(\begin{array}{c}
y_1 \\ y_2\\ y_3\\ y_4\\m_1\\m_2\\ m_3\\ m_4\\ m_1'\\ m_2'\\ m_3'\\ m_4'\\
\end{array}\right)
{\mbox{\large$=$}} \left(\begin{array}{c}
0 \\ 0\\ 0\\ 0\\ 0\\ 0\\ 0\\ 0\\0\\0\\ 0\\ 0\\
\end{array}\right)
\end{equation}
The solution set of this homogeneous linear equation system has dimension at
 least two because every vector $({\mathbf  y},{\mathbf  m}, {\mathbf  m}^\prime)$
 generated by  $(1,1,1,1,0,0,0,0, 0,0,0,0)$ and $(0,0,0,0,1,1,1,1,1,1,1,1)$ is
 a solution.
We claim that these are all the solutions.
This implies $y_1=y_2=y_3=y_4$
and $m_1=\dots=m_4'$, contradicting the fact that the lines $\pi_1,
\dots, \pi_4$ are distinct.

Let $M(r)$ be the $12\times 12$ coefficient matrix of (\ref{eq:12}).
Let us denote the characteristic polynomial of $M(r)$ by $f(r, \lambda)$.
We obtain
\begin{eqnarray*}
f(r,\lambda)&:=&\det(M(r)-\lambda I)\hfill \\
{} &=& {\lambda}^{12} -4 {\lambda}^{11} +{\lambda}^{10} (-r^2+5r)
+{\lambda}^9 (2 r^2-19 r+21)
+{\lambda}^8 (-6 r^2+40 r-68)\\
&{}&+\, {\lambda}^7 (r^4-3r^3+15 r^2-61 r+113) +{\lambda}^6 (2
r^4-30 r^3+107 r^2-99 r-94)
\\&{}&+\, {\lambda}^5 (2 r^5-26 r^4+120r^3-241 r^2+123 r+132)
\\&{}&+\, {\lambda}^4 (r^5-21 r^4+103 r^3-210 r^2+216 r-152)
\\&{}&+\, {\lambda}^3 (r^5-19 r^4+143 r^3-418 r^2+412r)
+{\lambda}^2 (6 r^4-38 r^3+60 r^2).
\end{eqnarray*}
Here the coefficient of $\lambda^2$ is not $0$ for $r\geq 4$
 (since $6 r^4-38 r^3+60 r^2=2r^2(r-3)(3r-10)$).
Thus the rank of $M(r)$ is $10$ and the solution set of
(\ref{eq:12}) has dimension 2, as stated.
The calculations have been verified by both the {\em Maple} and the
 {\em Mathematica} programs. \qed

\subsection{$3\times 3$ Euclidean grids}

When $r=3$ there are crossing families of straight lines.
Let $\mathbf  y$, $\mathbf  m$, ${\mathbf  m}^\prime\in \bR^3$, as before
 with $y_1> y_2>  y_3$, such that
 $\pi_i:=\{ (1,y_i), (2,y_i+m_i), (3,y_i+2m_i)\}$,
 $\rho_i:=\{ (j,y_i+(j-1)m_j^\prime): 1\leq j\leq 3\}$. Then
 $\cP:=\{ \pi_1, \pi_2, \pi_3\}$, $\cR:=\{ \rho_1, \rho_2, \rho_3\}$
 may form a crossing pair of families, i.e.,
$V_j:=\{\ell_j\cap \pi_i: 1\leq i\leq 3\}= \{\ell_j\cap \rho_i: 1\leq i\leq 3\}$
 with $\Pi:=\cup V_j$, $|V_1|=|V_2|=|V_3|=3$.
For example, for any $y,m,a,b$ (with $a,b > 0$) 
 we have such a system with values
\begin{multline}\label{eq:3}
 {\quad}{\mathbf  y}=(y+4a+2b, y-2a+2b, y-2a-4b). \\
  {\mathbf  m}=(m-3a, m-3b, m+3a+3b), \quad  {\mathbf  m}^\prime=(m-3a-3b, m+3a, m+3b).\quad{}
\end{multline}
The corresponding crossing systems are
\begin{eqnarray*}
\pi_1&=& \{ y+4a+2b , y+m+ a+2b, y+2m-2a+2b\}\\
\pi_2&=& \{ y-2a+2b , y+m-2a-b, y+2m-2a-4b\}\\
\pi_3&=& \{ y-2a-4b , y+m+ a-b,\enskip y+2m+4a+2b\}
  \end{eqnarray*}
  and
\begin{eqnarray*}
\rho_1&=& \{ y+4a+2b , y+m+ a-b,\enskip y+2m-2a-4b\}\\
\rho_2&=& \{ y-2a+2b , y+m+a+2b, y+2m+4a+2b\}\\
\rho_3&=& \{ y-2a-4b , y+m-2a-b, y+2m-2a+2b\}.
\end{eqnarray*}

\begin{lemma}\label{le:210}
If the set of slopes $M:=\{m_1, m_2, m_3\}\cup \{m_1', m_2', m_3' \}$ is
 $A_4$-free and $A_6$-free, then 
 $\cP$ and $\cR$ cannot be a set of $3\times 3$ crossing Euclidean lines.
 \end{lemma}

\noindent{\bf Proof.}\quad
From Lemma \ref{le:cross} we know that a
grid for $r=3$ has the following intersection pattern:

\[
\begin{array}{ccccc}
(1,y_{1}) &
\mbox{\begin{picture}(35,11)\thinlines\put(0,3){\line(1,0){35}}
\end{picture}}
& (2,y_{12}) &
\mbox{\begin{picture}(35,11)\thicklines\linethickness{2pt}\put(0,3){\line(1,0){35}}
\end{picture}}
& (3,y_{13}) \\
~ & \mbox{\begin{picture}(35,35)\thicklines\linethickness{2pt}
\put(0,0){\line(6,5){35}} \put(35,0){\line(-6,5){35}}\end{picture}}
& ~ & \mbox{\begin{picture}(35,35)\thinlines
\put(0,0){\line(6,5){35}} \put(35,0){\line(-6,5){35}}\end{picture}}
& ~ \\
(1,y_{2}) & ~ & (2,y_{22}) & ~ &(3, y_{23})\\
~ & \mbox{\begin{picture}(35,35)\thinlines \put(0,0){\line(6,5){35}}
\put(35,0){\line(-6,5){35}}\end{picture}} & ~ &
\mbox{\begin{picture}(35,35)\thicklines\linethickness{2pt}
\put(0,0){\line(6,5){35}}
\put(35,0){\line(-6,5){35}}\end{picture}}\\
(1,y_{3}) &
\mbox{\begin{picture}(35,11)\thicklines\linethickness{2pt}\put(0,3){\line(1,0){35}}
\end{picture}} & (2,y_{32}) &\mbox{\begin{picture}(35,11)\thinlines\put(0,3)
{\line(1,0){35}}
\end{picture}}  & (3,y_{33})
\end{array}
\]
Considering the intersection points
 $Q_{1,2}$, $Q_{1,3}$, $Q_{2,2}$, $Q_{2,3}$, $Q_{3,2}$, and $Q_{3,3}$,
 we derive the following system of linear equations.
\begin{equation}\label{eq:17}
\begin{array}{ccc ccc c ccc ccc}
y_1&&& +m_1&&& =&&y_2&&  &+m_2'&     \\
&&y_3& &&+2m_3& =&&y_2&&  &+2m_2'&  \\
&&y_3& &&+m_3&=&y_1&&&  +m_1'&&      \\
y_1&&& +2m_1&&&=&&&y_3&  &&+2m_3'   \\
&y_2&& &+m_2&&=&&&y_3&  &&+m_3'       \\
&y_2&& &+2m_2&&=&y_1&&&  +2m_1'&&   \\
\end{array}
\end{equation}

It is easy to see that all solutions of (\ref{eq:17})  are of the form of 
(\ref{eq:3}). \qed

\section{Constructions}

\subsection{Grid-free systems mod $q$}\label{ss:31}

In this section we prove our main results, Theorems~\ref{th:1} and~\ref{th:main}.
Given integers $q\geq r\geq 2$, $M \subset \{ 0,1,\dots, q-1\}$ define the hypergraph
 $\cF_M$ as follows.
Define the vertex set on the Euclidean plane
$$V:=\{(j,y):~1\le j\le r, ~y\in Z_q\}.$$
Thus $|V|=rq$, and let $V_j=\{(j,y):~y\in Z_q\}$.
For integers  $0\le y,m<q$ define the $r$-set
$$A(y,m)=\{(1,y), (2,y+m), \dots , (r,y+(r-1)m)\},$$
where the second coordinates are taken modulo $q$.
For $|m|<q/(4r)$  and $q/4< y< 3q/4$  or for $0\le y\leq y+(r-1)m<q$ the points of $A(y,m)$ are collinear.
For a subset of {\em slopes} $M\subseteq \{0,\dots ,q-1\}$ let
\begin{equation} \label{egyenes}
    {\cal F}_M:={\cF}_{M,q}^r=\{A(y,m):~y\in Z_q,~m\in M\}.
    \end{equation}
Obviously, this hypergraph is $r$-uniform, $|M|$-regular. Even  more, it can be decomposed into
 $|M|$ perfect matchings.

\begin{lemma} \label{le:31}
Suppose that  $q\geq r\geq 4$, and for all slopes
  $m\in M$ we have $|m|< q/(4r)$ (taken modulo $q$).
Then ${\cal F}_M$ is $\Grr$-free.
\end{lemma}

\noindent{\bf Proof.}\quad
Suppose  that given a grid ${\cal P}=\{A(y_i,m_i)\}_{i=1}^r, {\cal
R}=\{A(y_i^\prime ,m_i^\prime)\}_{i=1}^r\subset \cF$.
The hypergraph ${\cal F}_M$ is shift invariant, so we
may assume that $\exists \pi_1\in{\cal P}$ and $\exists \rho_1\in{\cal
R}$ which start in $(1,\lfloor q/2\rfloor)$.
In other words, replace each $A(y,m)$ by $A(y+\lfloor q/2\rfloor -y_1 ,m)$
 (the second coordinates are always taken modulo $q$).
The shifted system of hypergedges  have the same intersection structure as  $\cP$ and $\cR$.
So from now on, we may suppose that $y_1=\lfloor q/2\rfloor$.
Since the slopes $m_1$ and $m_1^\prime $ are small
 (i.e., $|m_1|, |m_1^\prime|< q/(4r)$) the points of
 $\pi_1:=A(\lfloor q/2\rfloor, m_1)$ and
  $\rho_1:=A(\lfloor q/2\rfloor, m_1^\prime)$ are forming
 Euclidean lines.
All other sets of $\cP$ and $\cR$ meet  either $\pi_1$ or $\rho_1$, and all other slopes
 are small (at most $q/(4r)$) so all other members of $\cP\cup\cR$ are forming Euclidean lines, too.

Finally, Lemma~\ref{le:28} completes the proof that $r$ should be at most $3$.
\qed

\begin{lemma}
If $q$ is a prime (and $q\geq r$), then ${\cal F}_M$ is a linear hypergraph.
\end{lemma}

\noindent
\noindent{\bf Proof.}\quad Well-known and easy.  $|A(y,m)\cap
A(y^\prime,m^\prime)|\ge 2$ implies that there exist $1\leq i\neq j\leq r$
\begin{eqnarray*}
y+im &\equiv & y^\prime+im^\prime \pmod{q}\\
y+jm &\equiv & y^\prime+jm^\prime \pmod{q},
\end{eqnarray*}
implying $(j-i)(m^\prime-m)\equiv 0 \pmod{q},$ a contradiction. \qed

\begin{lemma}\label{le:33}  If $q$ is a prime, $q> r^{4r}$, $r\geq 4$,  then the whole
   ${\cal F}_M$ with $M=Z_q$  is $\Grr$-free.
\end{lemma}

\noindent{\bf Proof.}\quad
Suppose  that the families  ${\cal P}=\{A(y_i,m_i)\}_{i=1}^r$ and $ {\cal
R}=\{A(y_i^\prime ,m_i^\prime)\}_{i=1}^r\subset \cF$ form a grid $\Grr$.
Apply  Minkowksi's theorem of simultaneous approximation (Lemma~\ref{le:mink})
 for the vector ${\mathbf  n}:=(m_1, \dots, m_r, m_1^\prime. \dots . m_r^\prime )\in \bR^{2r}$.
There exists an $0<\alpha <q$ such that for all $i\in [r]$ we have
   $|\alpha m_i|$ and $|\alpha m_i^\prime| \leq  q^{1-1/2r}< q/(4r)$ (mod $q$).
The collections $A(\alpha y_i,\alpha m_i)$ and
$A(\alpha y_i^\prime,\alpha m_i^\prime)$,
$i=1,\dots ,r$ have the same intersection pattern, i.e., they form a grid, too.
Then Lemma~\ref{le:31} implies that $r\leq 3$. \qed

\subsection{Grid-free systems for all $n$, the proof of Theorem~\ref{th:1}}\label{ss:32}

We use induction on $n$ to show that
  $\frac{n(n-1)}{r(r-1)}-c_rn^{8/5}<  \ex(n, \{ \I2, \Grr\})$ holds for $r\ge 4$ with
appropriate $c_r>0$.
Let $q$ be the largest prime $q\leq n/r$.
It is well-known~\cite{IP}  that $q> n/r -Cn^{3/5}$, where $C$ is an absolute constant.
Let $V_1, \dots, V_r$ be disjoint $q$-sets, and let $\cF_M$ be the grid-free
 hypergraph of size $q^2$ given by Lemma~\ref{le:33}.
Consider a grid-free, linear hypergraph $\cH$ on $q$ vertices.
By induction hypothesis, there is such a hypergraph of size  $|\cH|>\frac{q(q-1)}{r(r-1)}-c_rq^{8/5} $.
Put a copy of $\cH$, $\cH_j$, into each $V_j$ after randomly permuting their vertices,
 ($1\leq j\leq r$).
The union of these hypergraphs, $\cF :=\cF_M\cup \cH_1\cup \dots, \cup\cH_r$, is obviously
 linear.

Suppose that $\cP=\{ A_1, \dots, A_r\}$ and $\cR=\{ B_1, \dots, B_r\}$ are
 forming a grid in $\cF$.
Since $\cF_M$ is grid-free there should be an edge, say $A_j\in \cP$ such that
 $A_j\in \cH_j$.
The vertices of $A_j$ are covered by the edges of $\cR$, each meeting
 $A_j$ in a distinct singleton, so $\cR\subset \cF_M\cup \cH_j$.
If all $B_i\in \cF_M$ then, since $\cup \cP=\cup \cR$, we obtain that
 $A_i:=(\cup \cR)\cap V_i$ belongs to $\cH_i$ for each $i\in [r]$.
Call such a grid of {\em type 0} and their number is
 denoted by $g_0:=g_0(\cF)$.

Otherwise, there is an edge $B_j\in \cR\cap \cH_j$.
We claim that this edge $B_j$ is unique. If $B_j$ and $B_j'\in \cR\cap \cH_j$,
then all $A_i$ meet $V_j$ in at least two vertices, hence all $A_i\in \cH_j$,
hence $\cP$ and then $\cR$ are subfamilies of $\cH_j$, but $\cH_j$
 is grid-free.
It follows that $\cP\setminus \{ A_j\}$ and $\cR\setminus \{ B_j\}$
 are parts of $\cF_M$ and they are forming an $(r-1)\times (r-1)$ grid
 (on $V\setminus V_j$).
Call this $\cP\cup\cR$ a grid of {\em type $j$} and denote the
 number by type $j$ grids by $g_j:=g_j(\cF)$.

We obtain a grid-free family of size at least $|\cF|-g_0-g_1-\dots -g_r$
 if we leave out an edge from each grid in $\cF$.
Next we estimate the expected size of $g_j(\cF)$ when the vertices of each $\cH_i$
 are permuted randomly and independently.

Concerning type 0, there are at most ${|\cF|\choose r}$ choices of the
  vertex disjoint $B_1, \dots, B_r\in \cF_M$.
Given $B_1, \dots, B_r$ and  $j$  the probability that $A_j:=V_j\cap (\cup B_i)$
 indeed belongs to $\cH_j$ is exactly $|\cH_j|/\binom{q}{r}$,
 which is at most $1/\binom{q-2}{r-2}$.
These events are independent, and the number of ways to choose $B_1,\dots, B_r$
 is at most $\binom{q^2}{r}$, so the expected number of type 0 grids
$$
  \bbE(g_0) \leq  \binom{q^2}{r}  \binom{q-2}{r-2}^{-r} = O(q^{4r-r^2}).
  $$

Given $\cF_M$ one can count type $1$ grids as follows (the cases
$j>1$ are similar). Choose the edges $A_2$, $A_3\in \cF_M$. Each of
$B_2, \dots, B_r$ intersect both of them. Let $e_i$, $i=2,\dots ,r$
be the pair joining a vertex in $A_2\setminus V_1$ and $A_3\setminus
V_1$ if $B_i$ intersects them in these two vertices. Clearly, the
edges $e_i$ form a matching between the vertices $A_2\setminus V_1$
and $A_3\setminus V_1$ and they determine $B_2, \dots, B_r$.
Similarly, a matching between $B_2\setminus (V_1\cup A_2\cup A_3)$
and $B_3\setminus (V_1\cup A_2\cup A_3)$ determines $A_4, \dots,
A_r$. Finally, choosing  a vertex $c\in V_1$ the number of possible
choices so far is
$$
  \binom{|\cF_M|}{2}(r-1)!(r-3)!\times q =O(q^5).
  $$
The vertices $A_i\cap V_1$ and $B_i\cap V_1$  are called $a_i$ and $b_i$, resp.,
 $2\leq i\leq r$.
The probability that $\{c, a_2,\dots, a_r  \}$ and $\{c, b_2, \dots, b_r \}$
 both belong to $\cH_1$
 is at most
 $$ \frac{|\cH_1|}{\binom{q}{r}}\times \frac{(q-r)/(r-1)}{\binom{q-r}{r-1}}=O(q^{4-2r}),
 $$
  yielding $\bbE(g_1)=O(q^{9-2r})$.
Thus,  for $r\geq4$ the expected number of
 grids in $\cF$
 $$ \bbE \left(  g_0(\cF)+ \dots + g_r(\cF) \right)= O(q).
   $$
So there is a choice of permutations to make $\cF$ grid-free
 deleting only $O(q)$ edges.
This gives
$$
   \ex(n, \{ \I2, \Grr\})\geq |\cF|-O(q)\geq q^2+ r \frac{q(q-1)}{r(r-1)}-c_rq^{8/5}-O(q).
  $$
A short calculation shows that the right hand side is at least
$\frac{n(n-1)}{r(r-1)}-c_rn^{8/5}$ with some $c_r$. \qed

\subsection{Triangle-free systems, the proof of Theorem~\ref{th:main}}\label{ss:33}

Since the Tur\'an function is monotone we have to  consider only the case when
 $r$ divides $n$, $n=qr$.
Let $M\subset \{ 0,1, \dots, \lfloor q/(4r)\rfloor\}$ be an $r$-sum-free set of size
 $|M|> qe^{-\gamma _r\sqrt{\log q}}/(4r)$
provided by Lemma~\ref{le:behrend}
 and let $\cF_M$ be the family defined by~(\ref{egyenes}) in subsection~\ref{ss:31}.

By Lemma~\ref{le:31} for
$r\ge 4$ ${\cal F}_M$ is a linear hypergraph containing no grid.
Since the set of slopes $M$ is an $r$-sum-free set we cannot have three lines with
 slopes $m_1<m_3<m_2$ forming a triangle, either.
Otherwise, we get $c_1m_1+c_2m_2=(c_1+c_2)m_3$ for some $c_1+c_2\le r-1$
 (mod $q$, but it is not important here since all $|m_i|< q/(4r)$).
Finally, ${\cal F}_M=q|M|\ge n^2e^{-\beta_r\sqrt{\log n}}$
for some $\beta_r$, as stated.
 \qed

{\bf Proof} of Theorem~\ref{th:main} for $r=3$.\quad 
Here we establish the lower bound (\ref{eq:th_r=3}) 
  for $\ex(n,\{ \I2^3, \bbT_3, \bbG_{3\times 3} \})$
 when $n=3q$, and $q$ is a prime. 

First, recall that for $r=3$ ${\cal F}_M$ may contain a grid, see (\ref{eq:3}).
Therefore, to avoid grids and triangles at the same time we need
to choose the slopes in $M$ more restrictively.
For example, $M$ could be a set in $\{ 0, 1, \dots, \lfloor q/12\rfloor \}$
 which is $A_6$-free, $A_4$-free and $AP_3$-free simultaneously. 
Then  Lemma~\ref{le:210} gives that $\cF_M$ is linear and grid-free.
Also the $AP_3$-free property implies that ${\cal F}_M$ has no triangles either.
Lemma~\ref{le:A6} implies
$$
   \frac{2}{5}r_3(q/12)^{3/5}\times q \leq |M|q=|{\cal F}_M|\leq  
 \ex(n,\{ \I2^3, \bbT_3, \bbG_{3\times 3} \})
  $$
and then (\ref{eq:Sz}) completes the proof of the lower bound
 (\ref{eq:th_r=3}). \qed

Since ${\cal F}_M$ contains neither grids nor triangles it is
$3$-union free, too, by (\ref{eq:conj9}). 
This implies $n^{8/5-o(1)}< U_3(n,3)$.
Since  ${\cal F}_M$ is regular, (and can be split into matchings), it is an optimal 
 $3$-superimposed design.
The bound $|{\cal F}_M|/n=\Omega (n^{3/5-o(1)})$ exceeds
 the bound~(\ref{DyR})  for $k'(3,n)$ by D'yachkov and Rykov~\cite{DR2} for $r=3$.
Below we further improve both lower bounds with a different construction.

\subsection{Union-free triple systems}\label{ss:34}

Define the hypergraphs $\bbG_6$ and $\bbG_7$
 as follows on $6$ and $7$ vertices.\\ \indent
$E(\bbG_6):=\{ 123, 156, 426, 453\}$, \\ \indent
$E(\bbG_7):=\{ 123, 456, 726, 753\}$.\\
Note that both are three-partite and the 3-partition of their vertices is unique.

\begin{lemma}\label{le:3-partite_union} \quad
Suppose that $\cF$ is a three-partite, linear hypergraph. 
It is $3$-union-free if and only if it avoids $\G33$, $\bbG_6$, and $\bbG_7$.
  \end{lemma}

\noindent
{\em Proof:} \quad
We start like in the proof of Corollary~\ref{co:15}. 
Suppose, that ${\cal A}\neq{\cal B}$,  $|{\cal B}|\le |{\cal A}|\le 3 $,
 $\cup_{A\in {\cal A}}A=\cup_{B\in {\cal B}}B$ and $\bbG:=\cA\cup \cB$
 form a linear $3$-uniform hypergraph.
Then $\exists A_1\in{\cal A}\setminus {\cal B}$.
Since $|A_1\cap B|\leq 1$, to cover the elements of $A_1$ the family ${\cal B}$
 must contain $3$ sets. We obtain $|{\cal B}|=|{\cal A}|= 3$.
Moreover, the sets $B_1, B_2 , B_3\in \cB$ meet $A_1$ in distinct elements.

In the case of $\cA\cap \cB=\emptyset$ the latest property implies that every 
$a\in \cup \cA$ is covered by a unique $B\in \cB$, and every 
$b\in \cup \cB$ is covered by a unique $A\in \cA$, so $\bbG$ is 2-regular,
on 9 vertices, we obtain the grid $\G33$.

In the case of $|\cA\cap \cB|=2$, say $A_2=B_2$ and $A_3=B_3$ we have that 
$A_1\setminus (B_2\cup B_3)$ is a singleton and it must be the same element as
 $B_1\setminus (A_2\cup A_3)$. 
Taking into the account that $\bbG$ is 3-partite we obtain that 
 it is isomorphic to $\bbG_7$ when $A_2$ and $A_3$ are disjoint, 
 and it is isomorphic to $\bbG_6$ when $A_2$ and $A_3$ meet. 

Finally, in the case $|\cA\cap \cB|=1$, say $A_3=B_3$, the other four sets
 meet $A_3$ in exactly one element, so there is a vertex $v$ of $A_3$ 
 of degree at least $3$. 
If $v\in A_1\cap A_2\cap A_3$ then $B_1$ and $B_2$ covers $A_1\cup A_2\setminus v$
 and we could not finish because $\bbG $ is 3-partite.
Similarly, if $v\in A_2\cap B_2\cap A_3$ then $A_2\setminus v$ must be covered by $B_1$,
 a contradiction.  
So if the three configurations are avoided then $\cF$ is 3-union-free.  \qed

{\bf Proof of Proposition~\ref{pr:10}.}\quad 
The probabilistic lower bound 
 from Lemma \ref{le:erdos}
 and the previous Lemma imply the lower bound (\ref{eq:u3}) 
\begin{equation*}
\Omega(n^{5/3})\leq \ex(n, \{\G33, \bbG_6, \bbG_7\})\leq U_3(n,3). \qed 
  \end{equation*}

{\bf Proof of Proposition~\ref{pr:16}.}\quad 
Here we prove
 (\ref{eq:k3}) claiming $\Omega(n^{2/3})\leq k'(3,3n)$. 
The proof is a refined version of the previous proof.

Let $\cF$ be the three-partite complete hypergraph with parts 
 $V_1$, $V_2$ and $V_3$ where $V_i:=\{ (i,j): j\in Z_n\}$ 
as before.
Split it into $n^2$ perfect matchings 
$$ M(\alpha,\beta):=\{ \{(1,y), (2,y+\alpha), (3,y+\beta)\}: y\in Z_n \}
 $$
where the second coordinates are taken modulo $n$. 
We also call these {\em parallel classes}.

Choose independently each of the $n^2$ matchings with probability $p$, 
$p$ will be defined as $n^{-4/3}/2$. 
Call the obtained random hypergraph $\cH$.
Count the expected number of the arising 
 configurations $\I2^3$,  $\bbG_6$, $\bbG_7$, and $\G33$, those we want to avoid.
Here we have to be more careful, because although each edge
 belongs to $\cH$ with probability $p$ (so its expected size is $pn^3$)
 the choices of edges are not independent.
More precisely, the probability that a 
subhypergraph $\cA$ appears in $\cH$ is exactly $p^i$, where
 $i$ is the number of different parallel classes in $\cA$. 
 
Intersecting triples always belong to different classes, so they are independent,
so the expected number of $\I2$'s is $p^2\times 3n^2\binom{n}{2}$
 and the expected number of $\bbG_6$'s
 is  $p^4\times 2 \binom{n}{2}^3$.

The  expected number of grids in $\cH$ with independent hyperedges is $O(p^6n^9)$ 
 and the expected number of $\bbG_7$'s with independent edges is at most $O(p^4n^7)$.

A configuration $\bbG_7$ might have only one pair of dependent triples,
 namely the disjoint pair  $123$, $246$. The number of these
configurations is at most $3n^5$ (first we chose the triple corresponding to
$123$ in $n^3$ ways, then its parallel edge at most $n$ ways, and finally
 the 7th vertex at most $3n$ ways). So the expected number of 
 these is at most $3p^3n^5$. 
 
Consider, finally, the grids $\cA \cup \cB$ containing parallel triples.
Note that the parallel classes of $\cA$ and $\cB$ are distinct, because 
 every $A\in \cA$ meets every $B\in \cB$. 
The number of these configurations is $O(n^7)$
 (like before,  first we chose a triple in $n^3$ ways, then its parallel edge at most 
  $n$ ways, and finally  the three remaining  vertices at most $n^3$ ways).
So if the number of independent classes $i\geq 4$, then the expected number of 
 these grids is bounded by $O(p^4n^7)$.
 
If $i\leq 3$, then one of the three edges of the grid, say $\cA$, consists of 
 three parallel edges. The number of these configurations is at most $n^5$,
 so we have an upper bound $O(p^3n^5)$ for the expected number of these whenever $i=3$.
 
Finally, if $i=2$, then both $\cA$ and $\cB$ consists of parallel edges.
The number of those systems is at most $O(n^3)$ so we have an upper bound $O(p^2n^3)$.

Altogether the total expected number of configurations we want to avoid is
 bounded by a constant multiple of 
$$
 p^2n^4+ p^4n^6 + p^6n^9+p^4n^7+ p^3n^5+p^4n^7+p^3n^5+p^2n^3.
  $$
This is less than half of $E(|\cF|)=pn^3$ with our choice of $p$.
So we still have $\frac{1}{2}pn^3$ of the edges if we leave out 
 an edge from each of the configurations we want to avoid, and thus we make the
 rest 3-union-free.

Moreover, conveniently, we can erase the unwanted edges together with all its parallels
 (if $\cA$ is in $\cF$ then all of its shifted copy belongs to $\cF$), so 
 there is a random choice of $\cF$ where the remaining 3-union-free part
 is still large and union of matchings.  \qed

\subsection{Union free and cover-free graphs}\label{ss:35}

D'yachkov and Rykov \cite{DR2} showed  that there exists an optimal
$1$-superimposed code (i.e., a regular graph) with $\le {n^2}/{4}$ edges.
It is easy to see, that there are $k$-regular $n$-vertex graphs
 for all $k<n$ (for $nk$ even).
Indeed, if $n$ is even, and $n>k$ then, e.g., by Baranyai's
Theorem~\cite{B} the edge set of $K_n$ can be decomposed into
$n-1$ perfect matchings. Take $k$ of them. For odd $n$ (and $k$
even) $K_n$ can be decomposed into $(n-1)/2$
$2$-factors by a theorem of Tutte. Take $k/2$ of them.

D'yachkov and Rykov~\cite{DR2} showed that there exist optimal 2-superimposed
 designs (i.e., a $2$-union-free, $k$-regular graphs) for $k \le \log (n+2)-2$.
This can be improved to
\begin{equation}\label{eq:k2}
  k'(2,n)=\Theta(n^{1/2}).
  \end{equation}
Indeed, such an optimal design is just a $k$-regular, triangle and $C_4$-free
 graph.
Such a graph can be constructed using~(\ref{egyenes})
 by taking a mod $q$ Sidon set $M\subset Z_q$.
The defined $\cF_M$ is the desired bipartite graph.
We have an obvious upper bound
 $(n/2)k'(2,n)\leq \ex(n, \{ C_3, C_4\})$.
Since
$$
 \ex(n, \{ C_3, C_4\})\leq \ex(n, C_4)=(1+o(1))\frac{1}{2}n^{3/2}
$$
by~\cite{Brown, ERS} we have an upper bound $k'(2,n)\leq (1+o(1))\sqrt{n}$
 giving the right order of magnitude.
To determine the coefficient of the $\sqrt{n}$ seems to be a very difficult question.
Erd\H os and Simonovits~\cite{ES}  conjecture that
\begin{equation}
  (1+o(1))\frac{n^{3/2}}{2\sqrt 2}= \ex(n;\{C_4,C_3\}).\quad\quad (?)
  \end{equation}
They showed $\ex(n;\{C_4,C_5\})\sim n^{3/2}/(2\sqrt 2)$, i.e., forbidding  $C_5$ in magnitude is the
same as forbidding all non-bipartite graphs.

\section{Conclusion}

Our main result is that the widely investigated transversal design
$\cF_M$ (see~(\ref{egyenes})) with $M=Z_q$ is $\Grr$-free for $r\ge 4$. If $M$ is
$r$-sum-free then in addition $\cF_M$ has no triangles. It is
natural to ask what other small substructures can be avoided this
way.

\small

\bibliographystyle{amsalpha}

\begin{thebibliography}{AA}
\parskip=2pt plus 2pt
\lineskip=0pt plus 2pt

\bibitem{AZ} Nguyen Quang A, and T. Zeisel,
Bounds on constant weight binary superimposed codes,
 {\em Probl. of Control and Information Theory},
{\bf 17}, 1988, pp. 223--230.

\bibitem{AA1} N. Alon and V. Asodi, Tracing a single user, {\em European Journal of
Combinatorics}, {\bf 27} (8), 2006, pp. 1227--1234.

\bibitem{AA2} N. Alon and V. Asodi, Tracing many users with almost no rate penalty,
{\em IEEE Transactions on Information Theory}, {\bf 53} (1), 2007, pp.
437--439.

\bibitem{AS}
N. Alon and A. Shapira,
On an extremal hypergraph problem of Brown, Erd\H os and S\' os.
{\em Combinatorica} {\bf 26} (2006), no. 6, 627--645.

\bibitem{BS} L. Babai and V. T. S\'os, Sidon sets in groups and induced subgroups
of Cayley graphs, {\em Europ. J. Combinatorics,} {\bf 6}, 1985, pp. 101--114.

\bibitem{B} Zs. Baranyai, On the factorization of the complete uniform
hypergraph, {\em Proc. Colloq. Math. Soc. J\'anos B\'olyai},
Infinite and finite sets, Keszthely, Hungary, 1973.

\bibitem{Be} F. A. Behrend, On sets of integers which contain no three terms in
arihtmetical progression, {\em Proc. Nat. Acad. Sci. USA}, {\bf  32}, 331--333.

\bibitem{Bl}
A. Blokhuis, private communication, June 16, 2009.

\bibitem{DV} A. De Bonis and U. Vaccaro, Optimal algorithms for two group testing
problems and new bounds on generalized superimposed codes, {\em IEEE
Transactions on Information Theory,} {\bf 52} (10), 2006, pp. 4673--4680.

\bibitem{Br} A. E. Brouwer,
Steiner triple systems without forbidden subconfigurations,
{\em  Mathematisch Centrum Amsterdam}, ZW 104/77, 1977.

\bibitem{Brown}
   {W. G. Brown}, On graphs that do not contain a Thomsen graph,
   {\em Canad. Math. Bull.} {\bf 9} (1966), 281--289.

\bibitem{BES3}
W. G. Brown, P. Erd\H os, and V. T. S\'os,
On the existence of triangulated spheres in 3-graphs and related problems,
{\em Period. Math. Hungar.} {\bf  3} (1973) 221--228.

\bibitem{BESr}
W. G. Brown,  P.  Erd\H os, and V. T. S\'os,
Some extremal problems on $r$-graphs,
in: {\em New Directions in the Theory of the Graphs},
Proceedings of the Third Annual Arbor Conference on Graph Theory, Academic Press, New York, 1973, pp. 55--63.

\bibitem{CY}
Y. Caro and R. Yuster,
Packing graphs: the packing problem solved. 
{\em Electron. J. Combin.} {\bf 4} (1997),  Research Paper 1, 7 pp. (electronic). 

\bibitem{C}
C. J. Colbourn, private communication, Nov. 15, 2004.

\bibitem{CMRS}
C. J. Colbourn, E. Mendelsohn, A. Rosa, and  J. \v Sir\'a\v n,
Anti-mitre Steiner triple systems.
{\em Graphs Combin.} {\bf 10} (1994), no. 3, 215--224.

\bibitem{CR}
C. J. Colbourn and A. Rosa,
{\em Triple systems},
Clarendon Press, Oxford University Press, New York, 1999.

\bibitem{DR1} A. G. D'yachkov and V. V. Rykov, Bounds on the length of
disjunctive codes, {\em Problemy Peredaci Informacii,}
{\bf 18} (3), 1982, pp. 7--13.

\bibitem{DR2} A. G. D'yachkov and V. V. Rykov, Optimal superimposed codes
for R\'enyi's search model, {\em Journal of Statistical
Planning and Inference},  {\bf 100} (2), 2002, pp. 281--302.

\bibitem{E64}
 P. Erd\H os,
Extremal problems in graph theory, in: M. Fiedler (Ed.),
{\em Theory of Graphs and its Applications}, Academic Press, New York, 1964, pp. 29--36.

\bibitem{E64Isr}
 P. Erd\H os,
On extremal problems of graphs and generalized graphs.
{\em Israel J. Math.} {\bf 2} (1964), pp. 183–-190. 

\bibitem{ErdRoma}
P. Erd\H os, Problems and results in combinatorial analysis.
{\em in:} Colloquio Internazionale sulle Teorie Combinatorie (Rome, 1973), Tomo II, 3--17,
Atti dei Convegni Lincei, 17, Accad. Naz. Lincei, Rome, 1976.

\bibitem{EFF1} P. Erd\H os, P. Frankl, and Z. F\"uredi, Families of finite
sets in which no set is covered by the union of two others,
{\em Journal of Combinatorial Theory, Series A},  {\bf 33} (2), 1982, pp. 158--166.

\bibitem{EFF2} P. Erd\H os, P. Frankl, and Z. F\"uredi, Families of finite
sets in which no set is covered by the union of $r$ others,
{\em Israel J. of Mathematics},  {\bf 51} (1-2), 1985, pp. 79--89.

\bibitem{EFR}
P. Erd\H os, P. Frankl, and V. R\"odl,
The asymptotic number of graphs not containing a fixed subgraph and a problem
 for hypergraphs having no exponent,
{\em Graphs and Combin.}  {\bf 2} (1986) 113--121.

\bibitem{EF}
P. Erd\H os and Z. F\"uredi,
The greatest angle among $n$ points in the $d-$dimensional Euclidean space,
{\em Annals of Discrete Math.} {\bf 17} (1983),
{\it Combinatorial mathematics (Marseille-Luminy, 1981), North-Holland Math. Stud., 75}, pp. 275--283,
North-Holland, Amsterdam-New York, 1983.

\bibitem{ERS}
   {P. Erd\H os, A. R\'enyi, and V. T. S\'os},
    On a problem of graph theory,
   {\em Studia Sci. Math. Hungar.}
   {\bf 1} (1966), 215--235.

\bibitem{ES} P. Erd\H os and M. Simonovits,
Compactness results in extremal graph theory,
{\em Combinatorica}, {\bf 2}, 1982, pp. 275--288.

\bibitem{ET} P. Erd\H os and P. Tur\'an, On a problem of Sidon in additive
number theory and some related problems, {\em J. London Math. Soc.},
{\bf 16}, 1941, pp. 212--215.

\bibitem{EGy} T. Ericson and L. Gy\"orfi, Superimposed codes in
$R^n$, {\em IEEE Trans. Inform. Theory},  {\bf 34} (4), 1988, pp. 877--880.

\bibitem{FGG07}
A. D. Forbes, M. J. Grannell, and T. S. Griggs,
On 6-sparse Steiner triple systems,
{\em J. Combin. Theory Ser. A} {\bf 144} (2007) 235--252.

\bibitem{FGG09}
A. D. Forbes, M. J. Grannell, and T. S. Griggs,
Further 6-sparse Steiner triple systems,
{\em Graphs and Combin.} {\bf 25} (2009), no. 1, 49--64.


\bibitem{Fur25}
P. Frankl and Z. F\"uredi,
A new extremal property of Steiner triple systems,
{\em Discrete Mathematics\/}  {\bf 48} (1984), 205--212.

\bibitem{Fur51}
P. Frankl and Z. F\"uredi,
Union-free families of sets and equations over fields,
  {\em Journal of Number Theory\/} {\bf 23} (1986), 210--218.


\bibitem{Fur60} 
P. Frankl and Z. F\"uredi, 
Exact solution of some Tur\'an-type problems,
 {\em Journal of Combinatorial Theory, Ser.~A\/} {\bf 45} (1987), 226--262. 


\bibitem{Fuy}
Yuichiro Fujiwara,
Infinite classes of anti-mitre and 5-sparse Steiner triple systems,
{\em J. Combin. Des.} {\bf 14} (2006), no. 3, 237--250.

\bibitem{F} Z. F\"uredi, A note on $r$-cover-free families, {\em Journal of
Combinatorial Theory, Series A}, {\bf 73}, 1996, pp. 172--173.

\bibitem{FR} Z. F\"uredi and M. Ruszink\'o, Superimposed codes are
almost big distance ones, {\em Proc. 1997 IEEE
Int. Symp. Inform. Theory}, Ulm, Germany, June 29 - July 4,
1997, p. 118.

\bibitem{FuR} Z. F\"uredi and M. Ruszink\'o, An improved upper bound
of the rate of Euclidean superimposed codes, {\em IEEE Trans.
Inform.  Theory},  {\bf 45} (2), 1999, pp. 799--802.

\bibitem{GGW}
M. J. Grannell, T. S. Griggs, and C. A.  Whitehead,
The resolution of the anti-Pasch conjecture,
{\em J. Combin. Des.}  {\bf 8} (2000), no. 4, 300--309.

\bibitem{GM}
T. S. Griggs and J. P. Murphy, 101 Anti-Pasch Steiner triple systems of order 19,
{\em J. Combin. Math. Combin. Comput.} {\bf 13} (1993) 129--141.

\bibitem{GMP}
T. S. Griggs, J. Murphy, and J. S. Phelan, Anti-Pasch Steiner triple systems,
{\em J. Combin. Inf. Syst. Sci.} {\bf 15} (1990), 79--84.

\bibitem{HB} D. R. Heath-Brown, Integer sets containing no arithmetic progression,
{\em J. London Math. Soc.},  {\bf 35}, 1987, pp. 385--394.

\bibitem{H} F. K. Hwang, A method for detecting all defective members in
a population by group testing, {\em J. of the American Statistical
Association},  {\bf 67} (339), 1972, pp. 605--608.

\bibitem{HS} F. K. Hwang and V. T. S\'os, Non adaptive hypergeometric
group testing, {\em Studia Sci. Math. Hungar.}, {\bf 22}, 1987,
pp. 257--263.

\bibitem{IP}  H. Iwaniec and J. Pintz, Primes in short intervals. {\em Monatsh. Math.},
{\bf 98}, 1984, pp. 115--143.

\bibitem{J1} S. M. Johnson, On the upper bounds for unrestricted
binary error-correcting codes, {\em IEEE Trans. Inform. Theory},
{\bf 17} (4), 1971, pp. 466--478.

\bibitem{J2} S. M. Johnson, Improved asymptotic bounds for
error-correcting codes, {\em IEEE Trans.  Inform. Theory},
{\bf 9} (4), 1963, pp. 198--205.

\bibitem{J3} S. M. Johnson, A new upper bound for error-correcting codes,
{\em IRE Trans. Inform. Theory},  {\bf 8}, 1962, pp. 203--207.

\bibitem{KS} W. H. Kautz and R. C. Singleton, Nonrandom binary
superimposed codes, {\em IEEE Trans. Inform. Theory}, {\bf 10}, 1964,
pp. 363--377.

\bibitem{KL} H. K. Kim and V. Lebedev, On optimal superimposed codes,
{\em Journal of Combinatorial Designs}, {\bf 12} (2), 2004, pp. 79--91.

\bibitem{KLO} H. K. Kim, V. Lebedev, D. Y. Oh, Some new results on superimposed codes,
{\em Journal of Combinatorial Designs}, {\bf 13} (4), 2005, pp. 276--285.

\bibitem{LCGG}
A. C. H. Ling, C. J. Colbourn, M. J. Grannell, and T. S. Griggs,
Construction techniques for anti-Pasch Steiner triple systems,
{\em J. London Math. Soc.} (2) {\bf 61} (2000), no. 3, 641--657.

\bibitem{Ling}
A. C. H. Ling,
A direct product construction for $5$-sparse triple systems,
{\em J. Combin. Des.} {\bf 5} (1997), no. 6, 443--447.

\bibitem{Macu}
A. J. Macula,
A simple construction of $d$-disjunct matrices with certain constant weights,
{\em Discrete Math.} {\bf 162} (1996), no. 1-3, 311--312.

\bibitem{M} H. Minkowski, {\em Geometrie und Zahlen}, Leipzig und Berlin, 1896.

\bibitem{RWfelold} D. K. Ray-Chaudhuri and R. M. Wilson,
The existence of resolvable block designs,
{\em Survey of combinatorial theory}
(Proc. Internat. Sympos., Colorado State Univ., Fort Collins, Colo., 1971),
North-Holland, Amsterdam, 1973, pp. 361--375.

\bibitem{OB}  K. O'Bryant,
A complete annotated bibliography of work related to Sidon sequences,
{\em Electronic Journal of Combinatorics} {\bf 11} (2004).

\bibitem{Ru} M. Ruszink\'o, On the upper bound of the size of the
$r$-cover-free families, {\em Journal of Combinatorial Theory, Series A},
{\bf 66} (2), 1994, pp. 302--310.

\bibitem{ruzsa} I. Z. Ruzsa, Solving a linear equation in a set of integers I,
{\em Acta Arithmetica},  {\bf 65}, 1993, pp. 259--282

\bibitem{RSz} I. Z. Ruzsa and E. Szemer\'edi, Triple systems with no six points carrying
three triangles, in {\em Combinatorics, Keszthely, 1976,
Coll. Math. Soc. J. Bolyai}, {\bf  18} Volume II., pp. 939--945.

\bibitem{SS2}
G. N. S\'ark\"ozy and S. Selkow,
An extension of the Ruzsa-Szemer\'edi theorem.
{\em Combinatorica} {\bf 25} (2005), no. 1, 77--84.

\bibitem{SSk}
G. N. S\'ark\"ozy and S. Selkow,
On a Tur\'an-type hypergraph problem of Brown, Erd\H os and T. S\'os.
{\em Discrete Math.} {\bf 297} (2005), no. 1-3, 190--195.

\bibitem{S} V. T. S\'os,  An additive problem in different structures,
{\em Proc. of the Second Int. Conf. in Graph Theory,
Combinatorics, Algorithms, and Applications}, San Fra. Univ.,
California,  July 1989. SIAM, Philadelphia, 1991, pp. 486--510.

\bibitem{Sz} E. Szemer\'edi, Integer sets containing no arithmetic progression, {\em
Acta Math. Hungar.},  {\bf 56}, 1990, pp. 155--158.

\bibitem{T}
Luc Teirlinck, review MR2455590 (2009i:05037) on ~\cite{W08}
 {\em Mathematical Reviews} 2009.

\bibitem{Wilson} R. M. Wilson, On existence theory for pairwise balanced 
 designs, I, II, III, {\it J. Combinatorial Th., Ser. A} {\bf 13} (1972), 
 220--245, 246--273, {\bf 18} (1975),  71--79. 

\bibitem{W05}
A. Wolfe,
5-sparse Steiner triple systems of order $n$ exist for almost all admissible $n$,
{\em Electron. J. Combin.} {\bf 12} (2005), Research Paper 68, 42 pp.

\bibitem{W06}
A. Wolfe,
The resolution of the anti-mitre Steiner triple system conjecture.
{\em J. Combin. Des.} {\bf  14} (2006), no. 3, 229--236.

\bibitem{W08}
A. J. Wolfe,
The existence of 5-sparse Steiner triple systems of order $n\equiv 3\, (\bmod 6)$, $n\notin \{9,15\}$,
{\em J. Combin. Theory Ser. A} {\bf 115} (2008), no. 8, 1487--1503.

\end{thebibliography}

\end{document}